  \documentclass[12pt]{amsart}

  \usepackage[margin=1in]{geometry}

\usepackage{amssymb}
\usepackage{graphicx}
\usepackage{url}
\graphicspath{{images/}{../ApprDeepComps/images/}}
\usepackage[utf8]{inputenc}
\usepackage[english]{babel} 
\usepackage{amsmath}
\usepackage{amsfonts}
\usepackage{enumerate}
\usepackage{mathrsfs}
\usepackage{mathtools}
\usepackage[mathcal]{eucal}
\usepackage{comment}
\usepackage{color}
\usepackage{csquotes}
\usepackage{footmisc}

\counterwithin*{footnote}{section}

\numberwithin{equation}{section}

  \theoremstyle{plain}
  
  \newtheorem{therm}{Theorem}
  \newtheorem{proposition}[therm]{Proposition}

  \newtheorem{example}[therm]{Example}
  
  \theoremstyle{remark}
  \newtheorem{remark}[therm]{Remark}
  \newtheorem*{remark*}{Remark}
  \newtheorem{remarks}[therm]{Remarks}

\DeclareMathOperator{\T}{T\!}
\DeclareMathOperator{\tp}{tp}
\DeclareMathOperator{\Ilim}{\mathcal{I}lim}
\DeclareMathOperator{\Ulim}{\mathcal{U}lim}
\DeclareMathOperator{\Vlim}{\mathcal{V}lim}
\DeclareMathOperator{\id}{id}
\DeclareMathOperator{\Fix}{Fix}

\newcommand{\kR}{k\ensuremath{{}_{\RR}}}
\newcommand{\FFix}{\ensuremath{\mathtt{FindFix}}}
\newcommand{\ODEsolve}{\ensuremath{\mathtt{ODEsolve}}}

\newcommand{\Ftp}{\Phi\text{-}\!\tp}

\newcommand{\dd}{\ensuremath{\operatorname{d}}}
\newcommand{\dP}{\ensuremath{\dd_P}}
\newcommand{\ev}{\ensuremath{\operatorname{ev}}}

\newcommand{\nrm}[1]{\left|#1\right|}
\newcommand{\Nrm}[1]{\left\|#1\right\|}
\newcommand{\Ninf}[1]{\Nrm{#1}_{\infty}}

\newcommand{\CC}{\mathbb{C}}
\newcommand{\NN}{\mathbb{N}}
\newcommand{\RR}{\mathbb{R}}

\newcommand{\cC}{\mathcal{C}}
\newcommand{\cD}{\mathcal{D}}

\newcommand{\cI}{\mathcal{I}}
\newcommand{\cL}{\mathcal{L}}
\newcommand{\cP}{\mathcal{P}}
\newcommand{\cQ}{\mathcal{Q}}
\newcommand{\cT}{\mathcal{T}}
\newcommand{\cU}{\mathcal{U}}
\newcommand{\cV}{\mathcal{V}}

\newcommand{\Cp}{\mathrm{C}_{\mathrm{p}}}

\newcommand{\fD}{\mathfrak{D}}
\newcommand{\fL}{\mathfrak{L}}
\newcommand{\fT}{\mathfrak{T}}
\newcommand{\fX}{\mathfrak{X}}
\newcommand{\ff}{\mathfrak{f}}
\newcommand{\fg}{\mathfrak{g}}
\newcommand{\fid}{\mathfrak{i}}
\newcommand{\fu}{\mathfrak{u}}
\newcommand{\fv}{\mathfrak{v}}
\newcommand{\fx}{\mathfrak{x}}

\newcommand{\LSh}{\fL_{\mathrm{Sh}}}
\newcommand{\TShT}{T_{\mathrm{Sh}}}
\newcommand{\TSh}{\fT_{\mathrm{Sh}}}

\newcommand{\bv}{\mathbf{v}}
\newcommand{\bu}{\mathbf{u}}
\newcommand{\bw}{\mathbf{w}}

\newcommand{\rb}{r_{\bullet}}
\newcommand{\sbb}{s_{\bullet}}
\newcommand{\gammab}{\gamma_{\bullet}}

\begin{document}

\title{Approximability of deep computations}

\newcommand\PaperKeywords{Deep computations, ultracomputations, deep equilibrium
models, idempotent ultrafilters}

  \author{Samson Alva}
  \address{Department of Economics, The University of Texas at San Antonio}
  \email{samson.alva@utsa.edu}

  \author{Eduardo Dueñez}
  \address{Department of Mathematics, The University of Texas at San Antonio}
  \email{eduardo.duenez@utsa.edu}

  \author{José Iovino}
  \address{Department of Mathematics, The University of Texas at San Antonio}
  \email{jose.iovino@utsa.edu}

  \author{Claire Walton}
  \address{Department of Electrical Engineering and Department of Mathematics,
           The University of Texas at San Antonio}
  \email{claire.walton@utsa.edu}

  \date{\today}

\begin{abstract}
This is the first of a series of papers in which we study \emph{deep computations} (ultracomputations) and \emph{deep iterates}, formalizing the ideas of “asymptotic limit” of computations and compositional iterates, respectively.
In this first paper of the series, we characterize deep computations that are \emph{bona fide} computable, and prove the existence of \emph{deep equilibria}, which hitherto have been found only empirically in deep learning.
A subsequent paper will study the complexity of ultracomputations.
Our approach adapts and combines technology from the topology of function spaces, topological dynamics, and model theory. 
\end{abstract}


\maketitle

\section*{Introduction}
\label{sec:intro}

This is the first of a series of papers in which we combine methods of topology, measure theory, and model theory to study asymptotic behavior of computations as parameters of the computations tend toward a limit, e.g., the depth of a neural network tending to infinity, or the time interval between layers of a time-series network tending toward zero.
Recently, particular cases of this concept have attracted considerable attention in deep learning research (e.g., Neural ODEs/Ordinary Differential Equations~\cite{Chen-Rubanova-Bettencourt-Duvenaud:2018}, Physics-Informed Neural Networks~\cite{Raissi-Perdikaris-Karniadakis:2019}, and Deep Equilibrium/DEQ models~\cite{BaiKolterKoltun2020}).
The formal framework introduced here  provides a unified setting to study these limit phenomena from a foundational viewpoint.

While our methodology is abstract, it is motivated by the contemporary design features of deep (artificial) neural networks
—particularly, recurrent/recursive architectures where the parameters (weights and biases) of some layer(s) are reused to obtain multiple “weight-tied” or “cloned” new layers that are otherwise unique as such.
Both DEQ models and Neural ODEs replace explicit, layer-by-layer passes with asymptotic limits —either by driving a weight-tied transformation to an iterative stable fixed point in DEQ models, or by application of an unbounded number of “vanishingly small” residual layers (idealized as the continuous-time evolution under a differential equation) in the case of Neural ODEs.
From a computational perspective, the immense memory cost of back-propagation through an unbounded number of layers is effectively bypassed by delegating the forward and backward passes to generic “black-box” numerical algorithms (root-finding algorithms and ODE solvers, respectively).
Such methods bypass direct propagation by computing the limit by a different approach implicitly involving differentiability
—such considerations exceed our purely topological present approach, but the present framework may be extended in scope to capture differentiability.
An extended such framework is outlined in Appendix~\ref{sec:def-equi}.
At any rate, our study of continuity properties of asymptotic computations is a prerequisite to any analysis of their differentiability properties.
We find precise criteria characterizing “deep computation states” that are effectively computable
—a property we formalize as \emph{definability}.

The central notion of the paper is that of a \emph{deep computation}  (or \emph{ultracomputation}), that is, an accumulation point  of standard computations. Since deep computations are defined through nonconstructive processes (namely, ultrafilter limits), they are intrinsically  noncomputable; 
however, in engineering applications, they often happen to be effectively computable through standard libraries.
One of the two main results of the paper is a multipart characterization of  deep computations that are computable: 
We provide three different equivalent senses in which an asymptotic/deep computation may be regarded as effectively computable, at least in principle:
(\emph{i})~it is uniformly approximable by polynomials in quantities that are given \emph{a priori} as computable primitives
(i.e., it is definable in the sense of model theory);
(\emph{ii})~it is \emph{continuously extendable};
(\emph{iii})~certain iterated limits (intuitively, limits in state space and limits in time space) are exchangeable.
The proofs of these equivalences form the content of §~\ref{sec:DI-computability}, where the main result is Theorem~\ref{thm:fundl-thm-defin}.

The other main result of the paper is the formal proof of existence of \emph{deep equilibria}, which hitherto have been found only empirically in deep learning, by Bai--Kolter--Koltun \cite{BaiKolterKoltun2020}.
This material is in §§~\ref{S:deep-comp}–\ref{S:DEQ}, where the highlights are Proposition~\ref{prop:existence-ucomps} and Theorem~\ref{thm:DEQ-exist}.

We utilize and adapt technology from topology, combinatorics, and model theory. 
Since this is the first and foundational paper of a series, the results imported from these areas are classical:
The main topological tool used to characterize effectively computable deep computations is  Theorem~\ref{thm:Grothendieck}, which is an aggregation of several results of Grothendieck from the early 1950s~\cite{Grothendieck:1952};
the main combinatorial tool is Theorem~\ref{thm:DEQ-exist}, which is essentially a particular case of the classical Ellis--Numakura Lemma (see~\cite{Furstenberg:1981}); and the main ideas imported from model theory are the concept of \emph{space of types} (more precisely, under the view, pioneered by Krivine~\cite{Krivine:1976,Krivine-Maurey:1981}, of spaces of types as spaces of continuous real-valued functions, thus allowing us to connect spaces of types with topological function spaces), and crucially, Shelah's theory of stable theories~\cite{Shelah:1971,Shelah:1990}.
We attempt to make the paper broadly accessible and self-contained, assuming no particular technical expertise in those fields from the reader.%
\footnote{By contrast, our upcoming paper~\cite{DIMST}, which investigates the computational complexity of asymptotic computations, utilizes more recent and far deeper topological, combinatorial, model-theoretic and measure-theoretic tools:
the Todorčević trichotomy for Rosenthal compacta~\cite{Todorcevic:1999}, Shelah’s dependent and independent theories~\cite{Shelah:1971,Shelah:1990}, and recent classification results in topological dynamics and topological entropy obtained by Glasner and Megrelishvili~\cite{Glasner-Megrelishvili:2022}.}

To ease the introduction of the main ideas, in §~\ref{sec:extended-overview} we present a self-contained abridged summary, for the special case when our computation structures are countable, of the more general results expounded in §§~\ref{S:deep-comp}–\ref{sec:DI-computability}.

Readers with experience in model theory or in Banach space geometry will realize that the ideas presented here are strongly influenced by the work of C.~C.~Chang and H.~J.~Keisler on continuous model theory and model theory of real-valued structures~\cite{Chang-Keisler:1966, Keisler:2023} as well as the work of J.-L.~Krivine in Banach space theory~\cite{Krivine:1976,Krivine-Maurey:1981}.
We owe a great debt of gratitude to these giants
—we stand on their shoulders.

We are grateful to Frank Tall for his generous guidance through the world of $\Cp$-theory.


\section{Extended overview of the framework and main results}
\label{sec:extended-overview}

In this section we present a self-contained, abridged, and hopefully illuminating summary of the more general framework and results.
The material in §~\ref{sec:CCS-CSS} \emph{et seq.} does not otherwise depend on this overview section, which may be skipped entirely without logical harm.

For ease of exposition throughout this overview, we consider only computations whose states are characterized by (at most) countably many real-valued quantities (“features”).
All proofs presented in subsequent sections are omitted here.

\subsection{Definitions}
\label{sec:defs}

Fix a set $\cP = (P_n)_{n\in\NN}$ of countably many distinct symbols~$P_n$ called \emph{names of observables}, or \emph{distinguished predicate symbols}:
each symbol $P_n$ is a purely syntactic but distinct label.
(We could even take $\cP$ as the set~$\NN$ of natural numbers, regarding “$P_n$” simply as a syntactic formal alias for the element $n\in\NN$.)
Let $\RR^{\cP} = \prod_{Q\in\cP}\RR$ ($=\RR^{\NN}$) be the space of all functions $\cP\to\RR$, each regarded as a real tuple $v = (v_n)_{n\in\NN}$;
the space $\RR^{\cP}$ is endowed with the product topology, i.e., the topology of entry-wise convergence of such tuples.
Each $n\in\NN$ names a coordinate (projection) map $\pi_n(\cdot): \RR^{\cP}\to\RR: v \mapsto v_n$.
The real quantity $\pi_n(v) = v_n$ is called the \emph{$n$-th observable} of~$v$.

Fix an arbitrary nonempty subset $L\subseteq\RR^{\cP}$, which we shall call the \emph{state space};
its elements are called \emph{states}.
(We may also call these the \emph{layer state space} and \emph{layer states} to capture the neural-network intuition explained in the introduction.)
\emph{The (real-valued) observable on $L$ with symbol~$P_n$} is the map $P_n(\cdot) = \pi_n{\restriction}L: L\to\RR$ obtained by restricting $\pi_n$ to~$L$;
the \emph{(observable) $P_n$-feature} of $v\in L$ is the real value $P_n(v) = v_n$.
The pair $\underline{L} = \langle L, (P(\cdot))_{P\in\cP}\rangle$ is called a \emph{computation states structure (CSS)};
it will henceforth be denoted simply $\underline{L} = \langle L,\cP\rangle$ by an abuse of notation whereby we identify each symbol $P_n\in\cP$ with the corresponding feature~$P_n(\cdot): L\to\RR$.
(Such abuse of notation will be quite frequent throughout.)

The topological closure $\fL \coloneqq \overline{L}\subseteq\RR^{\cP}$ is called the \emph{space of (layer) state types of~$L$};
its elements are called \emph{state types}.
Elements $v\in L$ are called \emph{realized states} to distinguish them from state types $\bv\in\fL\setminus L$, called \emph{unrealized} (when such exist).

Each symbol $P_n\in\cP$ gives an observable feature (distinguished predicate)~$\fL\to\RR$ by restriction of the projection $\pi_n$;
it is the unique continuous extension of $P_n(\cdot): L\to\RR$ to~$\fL$, and will still be denoted $P_n(\cdot)$ (or even just~$P_n$).
The interpretations of $P_n$ on $L$ and on~$\fL$ are regarded as “the same observable” in practice, since either one determines the other.

The word “predicate” shall mean “real-valued function” and be used for functions on a topological space such as $L$ or~$\fL$;
however, such functions are otherwise completely arbitrary
—possibly discontinuous.
Thus, for instance, an arbitrary real function $f: L\to\RR$ is a predicate on~$L$, whereas observables are \emph{continuous} predicates (as immediate from the definition of the topology on~$L$).

A \emph{sizer} is a family $\rb = (r_P)_{P\in\cP}\in [0,\infty)^{\cP}$ of nonnegative reals indexed by observables $P\in\cP$.
Such a sizer names a compact subset $\RR{[\rb]}\coloneqq \prod_{P\in\cP}[-r_P,r_P] \subseteq \RR^{\cP}$.
Given a sizer~$\rb$, the \emph{$\rb$-shard of~$L$} (resp., \emph{of~$\fL$}) is $L[\rb] \coloneqq L\cap\RR{[\rb]}$ (resp., the closure $\fL[\rb]\coloneqq \overline{L[\rb]}\subseteq\RR^{\cP}$).
All type-shards $\fL[\rb]$ are compact (being closed in~$\RR{[\rb]}$).
Clearly, $\fL[\rb]\subseteq \fL\cap\RR{[\rb]}$ (equality need not hold).

\begin{proposition}\label{prop:kR-space}
  Let $\langle L,\cP\rangle$ be a CSS with countable observable collection~$\cP$.
  \begin{enumerate}
  \item the space $\fL$ of state types is metrizable;
  \item every state type is the limit of a sequence of realized states;
  \item every type $\bv\in\fL$ is shard-supported in the sense that $\bv\in\fL[\rb]$ for some sizer~$\rb$;
  thus, $\fL = \bigcup_{\rb}\fL[\rb]$ (where $\rb$ varies over all sizers);
  \item a predicate on~$\fL$ is continuous if its restrictions to arbitrary type-shards $\fL[\rb]\subseteq\fL$ are continuous.%
    \footnote{We thank F.~Tall for pointing out that $\fL$ is a \kR-space for $\cP$ countable.
      Indeed, property~(4) is a strengthening of the \kR\ property of~$\fL$ inasmuch as type-shards are compact
      (however, an arbitrary compact $K\subseteq\fL$ need not be included in any type-shard).}
\end{enumerate}
\end{proposition}
The proof uses metrizability in an essential way, hinting at the technical difficulties arising (from section~\ref{sec:CCS-CSS} on) when $\cP$ is possibly uncountable.
\begin{proof}
  The real line $\RR$ is topologized by the bounded metric $\dd(x,y)\coloneqq \rho(\nrm{y-x})$, where $\rho(t)\coloneqq t/(1+t)<1$.
  Since $\cP$ is countable, the space $\RR^{\cP}$ is metrizable, say by $\delta(\bu, \mathbf{v}) \coloneqq \sum_n2^{-n}\dd(P_n(u),P_n(v))<2$;
  therefore, its subspace $\fL$ is metrizable, proving~(1).
  By density of~$L$ in~$\fL$ and~(1), every type $\bv\in\fL$ is the limit of a sequence $v_{\bullet}=(v_n)\subseteq L$;
  hence, $K = v_{\bullet}\cup\{\bv\}$ ($= \overline{v_{\bullet}}$) is compact.
  The image $P(K)\subseteq\RR$ is bounded for each $P\in\cP$, hence $P(K)\subseteq[-r,r]$ for some $r = r_P>0$, and evidently $\fL[\rb]\supseteq K= \overline{v_{\bullet}} \ni \bv$ (where $\rb \coloneqq (r_P)_{P\in\cP}$).
  Assertions~(2) and~(3) follow.

  The compactness argument above is adapted to show that any convergent sequence $\bu_{\bullet} = (\bu_n)_{n\in\NN} \subseteq \fL$ is included in some type-shard.
  Indeed, for each $n\in\NN$, some sequence $v^{(n)}_{\bullet} \coloneqq (v^{(n)}_k)_{k\in\NN}\subseteq L$ satisfies $\lim_k v^{(n)}_k= \bu_n$;
  without loss of generality
  (upon replacing $v^{(n)}_{\bullet}$ by a sufficiently deep tail thereof if necessary),
  we may impose the following \emph{accelerated convergence} requirement: $\sup_{k\in\NN}\dd(v^{(n)}_k, \bu_n)\to 0$ as $n\to\infty$
  (the sequences $v^{(n)}_{\bullet}$ converge to their limits~$\bu_n$ “increasingly faster” as $n$ grows).
  Let $\bw \coloneqq \lim_n \bu_n$.
  Then, the set $K \coloneqq \{\bw, \bu_n, v^{(n)}_k: k,n\in\NN\} \subseteq \fL$ is compact:
  Given an open cover $\mathcal{G}$ of~$K$, we have $G\ni\bw$ for some $G\in \mathcal{G}$.
  By accelerated convergence, for all sufficiently large $n\in\NN$, say, for $n\ge N$, we have $G\supseteq\{\bu_n\}\cup v^{(n)}_{\bullet}$.
  For each $n<N$ there is also $G_n\in \mathcal{G}$ with $G_n\ni \bu_n$, hence $v^{(n)}_k\in G_n$ for all but finitely many $k\in\NN$;
  therefore, $\{G,G_n: n<N\} \subseteq \mathcal{G}$ covers all but finitely many points of~$K$, hence $\mathcal{G}$ has a finite subcover, so~$K$ is compact.
  Since each image $P(K)\subseteq \RR$ is (compact, hence) bounded, we deduce that $\fL[\rb] \supseteq K= \overline{\{v^{(n)}_k: n,k\in\NN\}} \supseteq \bu_{\bullet}$ for some sizer~$\rb$ (as before).

  Let now $\varphi:\fL\to\RR$ be discontinuous, say at $\bv\in\fL$;
  then, $\bv$ is the limit of some sequence $\bu_{\bullet} = (\bu_n)_{n\in\NN}$ in some type-shard $\fL[\rb]$ (by the preceding paragraph), but such that $\varphi(\bu_{\bullet}) \coloneqq (\varphi(\bu_n))\not\to  \varphi(\bv)$.
  We have $\bv\in\fL[\rb]$ (shards being closed in~$\fL$), so the restriction of~$\varphi$ to~$\fL[\rb]$ is discontinuous, proving~(4).
\end{proof}

A \emph{(syntactic) formula} is a purely formal real polynomial~$\varphi(P_1,\dots,P_k)$ in predicate symbols $P_1,\dots,P_k\in\cP$ (treated as pairwise commuting indeterminates).
Since each $P_i$ names a map $P_i(\cdot): \RR^{\cP}\to\RR$, such a formula~$\varphi$ itself names a \emph{polynomial function} (or just \emph{polynomial}) $\varphi(\cdot): \RR^{\cP}\to\RR$ called the \emph{interpretation} of~$\varphi$ which, in practice, we shall identify with the syntactic formula~$\varphi$.
(Different formulas may yield the same polynomial function, but this is not an issue in practice.)
By restriction of its interpretation on~$\RR^{\cP}$, a formula also gives polynomials on $L$ and on~$\fL$, as well as on shards $L[\rb]$ and $\fL[\rb]$;
moreover, by density of $L[\rb]$ in~$\fL[\rb]$, polynomials on $\fL[\rb]$ and $L[\rb]$ are in natural bijection
—and similarly for polynomials on~$\fL$ and $L$—
so we shall not distinguish between them.

A \emph{definable predicate} on a shard $L[\rb]$ is any predicate (i.e., real map) on $L[\rb]$ that is uniformly approximable by polynomials.
A definable predicate on $L$ is one whose restriction to every such shard $L[\rb]$ is definable.
Definable predicates on type space~$\fL$ and shard $\fL[\rb]$ are analogously defined.
Since a sizer $\rb$ specifies the same bounds on observable features of states $v\in L[\rb]$ and types $\bv\in\fL[\rb]$, and $L[\rb]$ is dense in~$\fL[\rb]$, one sees that definable predicates on $L[\rb]$ and $\fL[\rb]$ are identified (by the relations of restriction and extension, respectively), and similarly for definable predicates on $L$ and $\fL$.
Clearly, definable predicates on shards $L[\rb]$, $\fL[\rb]$ and on the space $L$ are continuous;
by Proposition~\ref{prop:kR-space}(4), definable predicates on $\fL$ are also continuous.

\label{transition,transform,t-t}
A map $L\to L$ (resp., $L\to \fL$, $\fL\to\fL$) is called a \emph{transition} (resp., a \emph{transform}, a \emph{transition-in-type (t-t)}).
By the inclusion $\fL\supseteq L$, every transition is a transform.
A transform~$f$ is \emph{extendable} if it extends to a continuous t-t on $\fL\supseteq L$.
A transform or t-t is  \emph{definable} if each of its features is definable.

For sizers $\rb,\sbb$, a transform (resp., a t-t) is called \emph{$\sbb$-confined on~$\rb$} if it restricts to a map $L[\rb]\to\fL[\sbb]$ (resp., $\fL[\rb]\to\fL[\sbb]$).
A set of transforms or t-ts \emph{is $\rb$-confined by $\sbb$} if each of its members is.

Any collection $\sbb^{[\cdot]} = (\sbb^{[\rb]})_{\rb}$ of sizers (itself indexed by sizers) is called a \emph{confiner}.
A transform or t-t is \emph{$\sbb^{[\cdot]}$-confined} if it is $\sbb^{[\rb]}$-confined on $\rb$ for all sizers~$\rb$.
Let $\cT, \fT$ be the sets of all confined transforms and t-ts, respectively;%
\footnote{More precisely, in this section $\cT$, $\fT$ denote the sets of \emph{confined} transforms and t-ts.
  In the systematic treatment (§~\ref{sec:t-t}), these same symbols shall denote the \emph{ambient} spaces $\cT = \LSh^{L}$ and $\fT = \LSh^{\LSh}$ of \emph{all} u-ts and t-ts whose confined subcollections are denoted $\cT[\sbb^{[\cdot]}]$, $\fT[\sbb^{[\cdot]}]$.}
the set of $\sbb^{[\cdot]}$-confined transforms (resp., t-ts) is denoted $\cT[\sbb^{[\cdot]}]\subseteq\cT$ (resp., $\fT[\sbb^{[\cdot]}]\subseteq\fT$).
Any subcollection of~$\cT[\sbb^{[\cdot]}]$ or $\fT[\sbb^{[\cdot]}]$ for some confiner $\sbb^{[\cdot]}$ is called \emph{uniformly confined} (or \emph{confined by $\sbb^{[\cdot]}$}).

A collection $R$ of sizers is called \emph{exhaustive} if $\fL = \bigcup_{\rb\in R}\fL[\rb]$.
A transform or t-t (resp., a set of transforms or t-ts) is called \emph{$R$-confined} if it is (resp., if all its members are) $\rb$-confined on~$\rb$, for all $\rb\in R$.
The set of $R$-confined t-ts is denoted $\fT[R]$.
(By exhaustiveness, $R$-confined transforms and t-ts are confined in the above sense.)

Any continuous t-t $\ff$ maps each shard $\fL[\rb]$ into some (compact subset of some) shard~$\fL[\sbb]$, so such $\ff$ is necessarily $\sbb^{[\cdot]}$-confined for some~$\sbb^{[\cdot]}$.

\begin{therm}\label{thm:extend-is-defin}
  A transform is extendable iff it is definable, in which case it is necessarily confined.
\end{therm}
\begin{proof}
  Let $f: L\to\fL$ be extended by a continuous $\ff: \fL\rightarrow\fL$.
  For fixed $Q\in\cP$, the “$Q$-th feature $Q{\circ}\ff$” of $\ff$ is continuous on the compactum~$\fL[\rb]$, hence uniformly approximable thereon by polynomials in observables $P_m(\cdot)$, by the Stone--Weierstrass Theorem
  (observables are continuous and separate points of~$\fL$);
  thus, $\ff$ is definable, and so is $f = \ff{\restriction}L$ \emph{a fortiori}.

  Conversely, if $f$ is definable, for fixed $Q\in\cP$, each restriction $Q{\circ}f{\restriction}L[\rb]$ of its $Q$-feature to an arbitrary state-shard $L[\rb]$ is a uniform limit of polynomials $\varphi$ in predicates.
 Each such $\varphi$ is (the restriction to~$L[\rb]$ of) a polynomial~$\boldsymbol{\varphi}$ on the compact type-shard~$\fL[\rb] = \overline{L[\rb]}$.
 Some sequence $(\boldsymbol{\varphi}_i)_{i\in\NN}$ of such polynomials converges uniformly on~$\fL[\rb]$ to a real function~$\mathbf{f}_Q^{[\rb]}$ on~$\fL[\rb]$ extending~$Q{\circ}f{\restriction}L[\rb]$ continuously.
 Since $\mathbf{f}_Q^{[\rb]}$ is continuous on the compactum~$\fL[\rb]$,
 it is bounded in magnitude thereon, say by $s = s_Q^{[\rb]}\in[0,\infty)$.
  Letting $Q$ vary, we obtain an $\sbb^{[\rb]}$-confined map $\ff^{[\rb]}: \fL[\rb]\to\fL[\sbb^{[\rb]}]: \bv\mapsto(\mathbf{f}^{[\rb]}_Q(\bv))_{Q\in\cP}$.
  Clearly, some (unique) $\sbb^{[\cdot]}$-confined $\mathbf{f}: \fL\to\fL$ extends all such~$\mathbf{f}^{[\rb]}$;
  such $\mathbf{f}$ is continuous since each entry $\mathbf{f}_Q$ is, by Proposition~\ref{prop:kR-space}(4).
\end{proof}

Theorem~\ref{thm:extend-is-defin} formalizes the (perhaps surprising) fact that \emph{non}-extendable transitions are not explicitly computable in terms of the observables~$P(\cdot)$.
  We remind the reader that the topology on~$L$ is the coarsest one for which all observables are continuous.
  However, even a \emph{continuous} transition~$f$, if non-extendable, is “uncomputable” in the sense that its features $P{\circ}f$ cannot be well approximated by continuous functions (e.g., polynomials) of observables~$Q(\cdot)$.
  Any sense of approximation \emph{cannot} be uniform;
  in fact, it cannot even be uniform on arbitrary shards~$L[\rb]$.

For that reason, extendability is a critical hypothesis in our main results.

\subsection{Deep computations and deep equilibria}
\label{sec:DC-DEQ}

Let $\sbb^{[\cdot]}$ be a confiner.
Recall that~$\cT[\sbb^{[\cdot]}]$, $\fT[\sbb^{[\cdot]}]$ are the sets of all transforms and t-ts, respectively, that are $\sbb^{[\cdot]}$-confined.
 (see page~\pageref{transition,transform,t-t} for the definitions).

An extendable transition will be called a \emph{computation}.

A \emph{compositional computation structure (CCS)} with countably many observables
\begin{equation*}
  \cC = \langle \underline{L}, \underline{\Gamma}, \ev\rangle
\end{equation*}
consists of:
\begin{itemize}
\item a CSS $\underline{L} = \langle L, \cP\rangle$ whose observables collection $\cP$ is countable;
\item a semigroup $\underline{\Gamma} = \langle\Gamma,\circ\rangle$, whose elements $\gamma\in\Gamma$ are called \emph{computations} of~$\cC$;
\item a continuous semigroup action $\ev: \Gamma\times L\to L$ of $\underline{\Gamma}$ on~$\underline{L}$.
\end{itemize}
Each computation $\gamma\in \Gamma$ gives a transition
\begin{align*}
  \gamma(\cdot): L &\to L\\
  v &\mapsto \gamma(v)\coloneqq\ev(\gamma,v).
\end{align*}
Under this identification, $\Gamma$ is a semigroup (under composition) of maps $L\to L$.

\emph{The CSS $\cC$ above will often be denoted simply $\cC = \langle \underline{L}, \underline{\Gamma}\rangle$ without explicitly naming the evaluation action~$\ev$, which, however, is always an implicit operation of $\cC$.%
\footnote{At any rate, the “functional application notation” $\gamma(x)$ for $\ev(\gamma,x)$ makes it essentially redundant to have a name for the action.}}

CCSs are required to satisfy the%
\footnote{When $\cP$ is uncountable, the Extendability Axiom takes a different form:
see~§~\ref{sec:Ax-Extend}.}\\
\emph{\textbf{Extendability Axiom.}\quad
  The transition $\gamma(\cdot)$ of any computation $\gamma\in \Gamma$ is extendable.}

By extendability, any computation is necessarily confined, so it may be regarded as a (confined) element $\gamma(\cdot)\in\cT$.

\subsubsection{Deep computations and ultracomputations}
\label{sec:DeepComp}

A \emph{deep computation (DC)} of a set $\Delta\subseteq\Gamma$ of computations is any confined transform $f\in\cT$ that is in the closure of (the set of transitions of) computations in~$\Delta$, in the topology of pointwise convergence.
Equivalently, a DC is any pointwise ultralimit $f \coloneqq \Ulim_i\gamma_i(\cdot)$ of any family $(\gamma_i)_{i\in I}\subseteq\Delta$ for any ultrafilter $\cU$ on the (otherwise arbitrary) index set~$I$, as long as each pointwise limit exists and the resulting map $f: L\to\fL$ is confined.
Although all transforms $\gamma(\cdot)$ of~$\Delta$ are deep computations, the name “deep computation” connotes a proper accumulation point that is otherwise not of the form~$\gamma(\cdot)$.

An \emph{ultracomputation (ucomp)} of~$\Delta$ is any (confined) accumulation point in~$\fT$ of the set $\tilde{\Delta}\subseteq\fT$ of transitions-in-type $\tilde{\gamma} = \gamma(\cdot)\in\fT$ extending computations $\gamma\in\Delta$.

In general,
\begin{itemize}
\item a DC need not be a map~$L\to L$
  —let alone need a ucomp restrict to such a map;
\item a DC need not have a unique extension to a ucomp.%
  \footnote{Clearly, any DC admits some extension to a ucomp, but such extension need not be continuous
   —nor, for that matter, be constructible in any explicit sense.}
\end{itemize}

\subsubsection{Deep iterates and deep equilibria}
\label{sec:DeepEq}

\paragraph{Deep iterates}
\label{sec:DI}

The (topological product) space $\fL^{\fL} = \prod_{\bv\in\fL}\fL$ of all t-ts $\ff: \fL\to\fL$ is a semigroup under composition.%
\footnote{Composition $(\ff,\fg)\mapsto\ff\circ\fg$ is continuous in the left argument~$\ff$, but generally \emph{not} in the right argument~$\fg$.
  In general, the set $\fT[\sbb^{[\cdot]}]$ of $\sbb^{[\cdot]}$-confined t-ts is not closed under composition.}
One sees that the subset $\fT\subseteq\fL^{\fL}$ is a sub-semigroup
(although its confined parts $\fT[\sbb^{[\cdot]}]$ are typically not closed under composition).

A \emph{deep iterate (DI)} of a computation $\gamma$ is any ultracomputation $\gamma^{(\cU)}\in\fT$ arising as ultralimit $\gamma^{(\cU)}\coloneqq \Ulim_n\tilde{\gamma}^{(n)}$ of iterates $\tilde{\gamma}^{(n)}\coloneqq \tilde{\gamma}\circ\dots\circ\tilde{\gamma}$ ($n$-fold) of the t-t $\tilde{\gamma}$ of~$\gamma$.
Note that the notion of deep iterate is strictly “in-type”, i.e., it is a transition-in-type —not a transition.
Being themselves confined by definition, deep iterates may be composed with any confined t-t.

\begin{proposition}[Cf., Propositions~\ref{prop:sh2sh-compactness}, \ref{prop:transf-type-spaces}, and~\ref{prop:existence-ucomps}]
  Fix any confiner~$\sbb^{[\cdot]}$, any exhaustive collection~$R$, and any set $\Delta\subseteq\Gamma$:
  \begin{enumerate}
  \item $\fT[\sbb^{[\cdot]}]$ is compact;
  \item $\fT[R]$ is a compact sub-semigroup of~$\fT$;
  \item Ultracomputations obtained from computations in~$\Delta$ form a closed sub-semigroup of~$\fT$;
  \item For any $\sbb^{[\cdot]}$-confined indexed family $(\gamma_i)_{i\in I}$, and any ultrafilter~$\cU$ on~$I$, the deep computation $\gamma_{\cU} \coloneqq \Ulim_i\gamma_i\in\cT[\sbb^{[\cdot]}]$ and the ultracomputation $\tilde{\gamma}_{\cU} \coloneqq \Ulim_i\tilde{\gamma}_i\in\fT[\sbb^{[\cdot]}]$ exist;
  \item If $\gamma$ is $R$-confined, then $\gamma$ has deep iterates, necessarily $R$-confined, of the form $\tilde{\gamma}^{(\cU)} = \Ulim_n\tilde{\gamma}^{(n)}$ for arbitrary nonprincipal $\cU\in\beta\NN$.
  \end{enumerate}
\end{proposition}

\paragraph{Deep equilibria}
\label{sec:DEQ}
A \emph{deep equilibrium (DEQ)} of a computation~$\gamma$ is an idempotent deep iterate~$\fid = \tilde{\gamma}^{(\cU)}$, i.e., one such that $\fid\circ\fid = \fid$.

\begin{therm}[Cf., Theorem~\ref{thm:DEQ-exist}]
  Let $R$ be any exhaustive collection.
  If $\gamma$ is an $R$-confined computation, then $\gamma$ has deep equilibria.
  In fact, one such DEQ is obtained as the ultralimit $\tilde{\gamma}^{(\cI)} = \Ilim_n\tilde{\gamma}^{(n)}$ from an arbitrary idempotent ultrafilter $\cI$ on~$\NN$.
\end{therm}

\subsection{Definability Criteria}
\label{sec:Grothendieck1}

  Ultracomputations, deep iterates and deep equilibria are typically \emph{not} definable, i.e., \emph{not} effectively computable
 —even when $\cP$ consists of a \emph{single} observable~$P$, let alone countably many!
  (Cf., Example~\ref{sec:unit-interval}~\emph{et seqq.})

  Theorem~\ref{thm:extend-is-defin} implies very strong restrictions on the ability to realize deep computations in any explicit fashion.
  One may ask for criteria ensuring that deep computations (or deep iterates, or deep equilibria) are effectively computable —i.e., definable.

\begin{therm}[Cf., Theorem~\ref{thm:fundl-thm-defin}]
\label{thm:definability}
  Let $\sbb^{[\cdot]}$ be a confiner, and let $\Delta$ be any collection of $\sbb^{[\cdot]}$-confined computations on a CCS with countable observables collection~$\cP$.
    Then, the properties below are equivalent:

    \begin{itemize}
    \item[(DD)]\emph{Deep Definability:}\quad
All deep computations of~$\Delta$ are definable (hence extend to continuous ultracomputations)---equivalently, conditions~(uExt) and~(UA) of Theorem~\ref{thm:fundl-thm-defin} hold.

\item[(LE)]\emph{Limit Exchange:}\quad
  For all predicates $P\in\cP$, all sizers $\rb$, and all sequences $v_{\bullet}\subseteq L[\rb]$ and $\gammab\subseteq \Delta$, the Limit Exchange identity
\begin{equation}\label{eq:limit-exchange1}
  \lim_m\lim_n P{\circ}\gamma_m(v_n) = 
  \lim_n\lim_m P{\circ}\gamma_m(v_n),
\end{equation}
    holds whenever the iterated limits on the left- and right-hand sides both exist.
    Moreover, in such case, each ultracomputation (hence, each DC) of~$\Delta$ is the (pointwise) limit of some sequence~$(\gamma_n)\subseteq\Delta$.
    The limit is attained uniformly on type-shards
    (\emph{a fortiori}, uniformly on state-shards).
  \end{itemize}
\end{therm}

\section{Structures for real-valued computations}
\label{sec:CCS-CSS}

This section, otherwise not dependent on the overview~§~\ref{sec:extended-overview}, begins the fully general and systematic development of the framework outlined therein.

We begin by introducing the notions of \emph{computation states structure (CSS)} and \emph{compositional computation structure (CCS)}, which lie at the foundation of our approach to real-valued computing.
The definitions of CSS and CCS in~§~\ref{sec:CSS} and~§~\ref{sec:CCS} are delicate but not difficult.
The preliminary informal discussion in §~\ref{sec:informal-structures} should demystify some of the formalism, but it is not an essential prerequisite to the material in §~\ref{sec:CSS} \emph{et seq.}

\subsection{Computations, states, observable features and predicates: A meteorological allegory}
\label{sec:informal-structures}
Consider physical quantities (such as temperature and barometric pressure) that are real-valued, and each of which may be observed at any given point.
For definiteness, consider points on or above the surface of earth, regarded as an idealized sphere.
A state $v$ captures the properties of a specific such point at a specific moment in time.
In such idealization, each physical quantity at any~$v$ is called a \emph{feature} of~$v$ (or \emph{observable feature} for emphasis).
Each such feature must be given a name (e.g., temperature, pressure, latitude, longitude, height, etc.);
these names are essential, for otherwise the real value of a feature of~$v$ is devoid of context.
We use the term \emph{observable} to refer to the \emph{name} given to any such property that may be observed;
in a formal treatment, we use (purely syntactic) symbols (e.g., “\texttt{T}” for temperature, “\texttt{p}” for pressure, “\texttt{lat}” for latitude, “\texttt{long}” for longitude, “\texttt{h}” for height, etc.)\ as observables.
An observable \emph{feature} of~$v$ is the value at~$v$ of the observable;
e.g., $v$ may have features $\mathtt{lat}(v) = +29.42$ (the \texttt{lat}-feature —i.e., latitude—of~$v$ is $29.42^{\circ}$~N), $\mathtt{long}(v) = -98.49$ ($v$ has \texttt{long}-feature —i.e., longitude— $98.49^{\circ}$~W),  $\mathtt{h}(v) = 229$ ($v$ is at $229$~m height), $\mathtt{T}(v) = 33.5$ (the temperature at~$v$ is $33.5^{\circ}$), etc.

We fix a symbol for each observable;
such symbols \texttt{P}, \texttt{Q}, \dots (not necessarily finitely many, or even countably many for that matter) will be called \emph{observable names} or \emph{distinguished predicate symbols}.
The set of distinguished predicate symbols (i.e., of symbols naming observables) will be denoted~$\cP$.
We shall denote the set of all possible states~$v$ by~$L$.
In the present discussion, $L$ might be taken to consist of points on the surface of our idealized spherical earth; it is perhaps more fitting to allow states $v\in L$ to refer to spatial points each at a specific moment in time.
Note that time $\mathtt{t}$ is \emph{not} an observable if one takes $L$ simply as the set of points on the sphere, but $\mathtt{t}$ is a valid observable on the set $L$ of states $v$ simultaneously encoding both location and time (in addition to other observables: temperature, pressure, etc.).

Any real-valued function on~$L$ is called a \emph{predicate}.
Each symbol $\mathtt{P}\in\cP$, at any state $v\in L$, has an associated real value~$P(v)$
(the switch to italic $P$ from typewriter-style $\mathtt{P}$ is a reminder that the symbol \texttt{P} has been “interpreted” to yield the actual value $P(v)$ of the \texttt{P}-\emph{feature} of~$v$).
Thus, the symbol~$\mathtt{P}$ entails a predicate (\texttt{P}-feature)
\begin{align*}
  P(\cdot): L &\to\RR\\
  v &\mapsto P(v),
\end{align*}
where the notation $P(\cdot)$ emphasizes the passage from the symbol~\texttt{P} to its interpretation.

Now that the distinction between an observable~\texttt{P} and the \texttt{P}-feature (distinguished predicate) $P(\cdot)$ interpreting it is clear, we shall henceforth use italic $P, Q, \dots$ simultaneously as formal (predicate) symbols denoting observables, and to denote the corresponding features;
in cases of potential confusion, we use the preferred notation $P(\cdot), Q(\cdot), \dots$ for predicates.
(Whenever $\cP$ is used as an index set, its members are regarded as symbols, never as predicates.)

Taking $L$ together with the features $P(\cdot)$ interpreting observables $P\in\cP$, we obtain a pair $\underline{L} \coloneqq \langle L, (P(\cdot))_{\mathtt{P}\in\cP}\rangle$ called a Computation States Structure (CSS) in~\ref{sec:state-types} below.
(In~$\underline{L}$, the collection of predicates $P(\cdot): L\to\RR$ is a family indexed by symbols~\texttt{P}.)
By an abuse of notation, we may denote such structure in the form $\underline{L} = \langle L,\cP\rangle$ wherein the collection $(P(\cdot):P\in\cP)$ of predicates is implicitly identified with the indexing set~$\cP$.

In the allegory, such features include the quantities \texttt{lat}$(v)$, \texttt{long}$(v)$ and \texttt{h}$(v)$, which are coordinates in the usual sense, as well as other features \texttt{T}$(v)$, \texttt{p}$(v)$ and time \texttt{t}$(v)$, which are not;
however, this suggests regarding the collection of \emph{all} features $(P(v))_{P\in\cP}$ of states as coordinatizing states~$v\in L$.
Each $P$-feature $P(v)$ is the “$P$-th coordinate” of~$v$ in an abstract sense;
the collection $\tp(v) \coloneqq (P(v))_{P\in \cP}$ is called the \emph{type of~$v$}.
Any state $v$ is uniquely characterized by its type.
A critical feature of our approach is to endow the state space~$L$ with the topology of “pointwise convergence”, i.e., a filter on (or: a sequence or net of) states converges to a state $v\in L$ iff the filter (or sequence, or net) of real-valued $P$-features converges to~$P(v)$, for each $P\in\cP$.

For the remainder of this subsection, we assume that \emph{the state space~$L$ is compact}.
In our allegory wherein height (and time) are observables allowed to take arbitrarily large values, compactness fails.
On the other hand, if the height and time intervals are \emph{bounded}
(e.g., $0 \le \mathtt{h}(v),\mathtt{t}(v) \le C$ for any fixed $C>0$),
the corresponding state space is compact.

On first approximation, a computation is a map $\gamma: L\to L$ transforming any given input state~$v$ to some output state~$\gamma(v)$.
(For simplicity, we use the same space~$L$ of input and output states.)
In our allegory, one may “visualize” computations as moving $v$ to another point~$\gamma(v)$, possibly at a different moment in time.
Maps $\gamma: L\to L$ should be considered “computable” in any reasonably explicit sense (say, by algorithms relying on arbitrary-precision arithmetic) only if output features $Q(\gamma(v))$ vary continuously with input features~$P(v)$, i.e., only when $\gamma$ is a continuous map~$L\to L$ in the topology of pointwise convergence of individual observable features.
Such requirement is consistent with the physics implied by our allegory.
We always require computations to be continuous.%
\footnote{Computations on a noncompact state space~$L$ are required to be \emph{extendable} in the sense of~§~\ref{sec:Sh-extend-comps}
  —a technical requirement significantly stronger than continuity.}

For illustration purposes, consider the “advance-time-by-1” computation $\alpha$ taking any state~$v$ of some point at some time \texttt{t}$(v) = t$ to the unique state~$w = \alpha(v)$ of the same point at time \texttt{t}$(w) = t+1$.
Features $T{\circ}\alpha, p{\circ}\alpha$ of the computation $\alpha$ give the temperature and pressure at a future moment $t+1$ in time from the state at present time~$t$.

When the state space~$L$ is compact, continuous computations~$\gamma: L\to L$ are effectively computable in a rather strong sense:
they are \emph{polynomially definable}.
This means that, up to any small fixed (but otherwise arbitrary) degree of precision, every output feature $Q(\gamma(v))$ is given (up to an error not exceeding the precision) by a polynomial on some input features~$P(v)$.
On the other hand (with apologies to meteorologists), our methods offer no insight on the specific polynomial approximating any output feature;
at any rate, such polynomial approximation would require fixing (finite) bounds for time~\texttt{t}.

As a by-product of choosing a common state space~$L$ both for computation inputs and outputs, computations are necessarily composable, i.e., any given computations naturally generate a semigroup of computations.
This gives rise to the notion of \emph{compositional computation structure (CCS)}, which is of the form
\begin{equation*}
  \cC = \langle \underline{L}, \underline{\Gamma}, \ev\rangle,
\end{equation*}
where $\underline{L}$ is a CSS, and $\underline{\Gamma} = \langle\Gamma,\circ\rangle$ is any semigroup under an (associative) composition operation~$\circ$, with elements $\gamma\in\Gamma$ representing computations on~$L$ via an \emph{evaluation} map $\ev: \Gamma\times L\to L: (\gamma,v)\mapsto\ev(\gamma,v)$ ($=\gamma(v)$, if $\Gamma$ is already a set of maps $\gamma: L\to L$).
Layer state transitions $\gamma(\cdot): v\mapsto \ev(\gamma,v)$ are assumed continuous on~$L$
(when $L$ is noncompact, we require them to be \emph{extendable} in the sense of~§~\ref{sec:Sh-extend-comps}).
CCSs are the natural structures to study compositions $\gamma_n\circ\gamma_{n-1}\circ\dots\circ\gamma_2\circ\gamma_1$ of $n$-many computations leading, as $n\to\infty$, to “deep computation states”, as well as “deep iterates” asymptotically approximated by $n$-fold iterates $\gamma^{(n)} = \gamma\circ\gamma\circ\dots\circ\gamma$ of a fixed computation~$\gamma$.

With suitable changes in definitions, our results apply to non-compact CSSs/CCSs.

\subsection{Computation States Structures}
\label{sec:CSS}
Fix an arbitrary nonempty set~$\cP$ whose members $P,Q,\dots$ will be called \emph{distinguished predicate symbols} or \emph{observable names}.

A \emph{Computation States Structure (CSS)} with observables~$\cP$ is of the form
\begin{equation*}
  \underline{L} = \langle L, (P(\cdot))_{P\in\cP}\rangle,
\end{equation*}
where
\begin{itemize}
\item $L$ is a nonempty set, called the \emph{sort (or space) of layer states};
\item For each symbol $P\in\cP$, the \emph{(distinguished) $P$-predicate%
      \footnote{In model theory of real-valued structures, an arbitrary real function is called a \emph{predicate}.}
      (or $P$-observable) of~$\underline{L}$} is a real function $P(\cdot): L\to\RR$.
\end{itemize}
We typically identify a symbol $P\in\cP$ with the predicate~$P(\cdot)$;
this entails a further abuse whereby we identify the predicate collection $(P(\cdot):P\in\cP)$ with $\cP$ itself;
thereby, the CSS above takes the form $\underline{L}=\langle L,\cP\rangle$.

\subsubsection{Types of states}
\label{sec:state-types}

In a CSS $\underline{L}=\langle L,\cP\rangle$, the \emph{type of a state~$v\in L$} is the indexed family $\tp(v) \coloneqq (P(v): P\in\cP)$ of its observable features.
Such type is called \emph{realized by~$v$};
it is a “vector” $\fv = (\fv_P)_{P\in\cP}$ with real entries $\fv_P = P(v)$ indexed by symbols~$P$.
Thus, such state types~$\fv$ are elements of the product (vector space)~$\RR^{\cP} = \prod_{P\in \cP}\RR$, which will always be regarded as topological product of copies of the real line~$\RR$ (endowed with its usual topology), one such copy for each~$P\in\cP$.
The topological subspace of realized types will be denoted $\tp(L) \coloneqq \{\tp(v): v\in L\}\subseteq\RR^{\cP}$.
(On the other hand, the \emph{linear} operations on the vector space $\RR^{\cP}$ will not play a direct role except in informal discussions —and in the Appendix.)

Ultrafilters on an infinite (“index”) set~$I$ will be denoted $\cU,\cV,\dots$;
we consider nonprincipal ultrafilters tacitly.
Given an ultrafilter $\cU$ on $I$, we say that an indexed family $\fv^{(\bullet)} \coloneqq (\fv^{(i)})_{i\in I}$ of elements $\fv^{(i)} = (\fv^{(i)}_P)_{P\in\cP}\in\RR^{\cP}$ \emph{converges to} $\fu=(\fu_P)_{P\in\cP} \in \RR^{\cP}$, or that $\fu$ \emph{is the $\cU$-ultralimit} of $\fv^{(\bullet)}$ if $\Ulim_i\fv_P^{(i)} = \fu_P$ for each $P\in\cP$
(i.e., when $\fu$ is the $\cU$-ultralimit of $(\fv^{(i)})$ in the  pointwise convergence topology—\emph{not} necessarily uniformly as $P$ varies).
Any such ultralimit $\fu$ is denoted $\Ulim_i\fv^{(i)}$ (it is unique if it exists).%
\footnote{When it exists, the $\cU$-ultralimit $s = \Ulim_i r_i$ of $(r_i)_{i\in I} \subseteq \RR$ is uniquely characterized by the following property:
  for every $\varepsilon>0$, the set $\{i\in I: \nrm{r_i-s} < \varepsilon\}$ belongs to~$\cU$ (i.e., is a “$\cU$-large” set).
  Ultralimits $\Ulim_iP(u^{{(i)}})$ may not exist since $\RR$ is not compact
  —however, all such ultralimits do exist when the values $P(u^{(i)})$ are bounded as $i$ varies.}

Elements~$\fu\in\RR^{\cP}$ arising as entry-wise ultralimits \emph{of realized types $\tp(v)$} in the above fashion (with $I$ and $\cU$ allowed to vary) are called \emph{types of (layer) states}, or \emph{ultrastates}.
Any realized state type is an ultrastate, but the converse fails in general.
The set $\fL$ of ultrastates is a closed subset $\fL = \overline{\tp(L)}$ (the bar denoting topological closure) of $\RR^{\cP}$, called the \emph{state type space}, and henceforth endowed with the subspace topology.
Since ultrastates need not be realized, the inclusion $\tp(L)\subseteq \fL$ is generally proper.

We shall adopt the convenient alternate notation $P(\fv)$ for the “$P$-th feature” $\fv_P$ of a type~$\fv\in\fL$, which treats $\fv$ as if it were realized
(i.e., as though $\fv$ were a state in~$L$).

\subsubsection{Topology on the layer state space}
\label{sec:L-topology}

We adopt a structural perspective wherein states are to be distinguished \emph{only} through their features;
thus, a state~$v\in L$ is implicitly identified with its type~$\tp(v)\in\RR^{\cP}$.
We topologize $L$ with (the “pullback” of) the product topology under such identification.
A slightly more concrete description of this topology is as follows:
For each predicate~$P\in\cP$, endow $L$ with the pseudometric $\dP(v,w) \coloneqq \nrm{P(w)-P(v)}$, and topologize $L$ by the collection $(\dP)_{P\in\cP}$ of all such pseudometrics.
This topology “by type” is the only one we shall introduce on~$L$
(except in certain examples meant to compare this topology to others, typically finer).

CSSs are assumed to satisfy the following:
\begin{itemize}
\item \emph{\textbf{Reduction Axiom for Computation States Structures.}
    States $v,w\in L$ are equal if their types $\tp(v), \tp(w)$ are equal.}
\end{itemize}

The Reduction Axiom above is equivalent to the requirement that distinct states be topologically distinguishable;
since $\RR^{\cP}$ is Hausdorff, reduction amounts to requiring that $L$ itself be a Hausdorff topological space (any two states have disjoint neighborhoods).

Even if not imposed \emph{a priori} on a CSS~$\underline{L}$, the Reduction Axiom is always satisfied if one replaces $L$ by its quotient $\tilde{L}\coloneqq L/{\tp}$ upon identifying equal-in-type states, and each predicate $P(\cdot)$ on~$L$ by the naturally induced predicate $\tilde{P}(\cdot)$ on~$\tilde{L}$.
From a structural viewpoint, $\langle L,(P(\cdot))_{P\in\cP}\rangle$ and $\langle \tilde{L}, (\tilde{P}(\cdot))_{P\in\cP}\rangle$ are identical
(isomorphic, in the sense of Keisler’s General Real-Valued Structures~\cite{Keisler:2023}).

\begin{remark}\label{rem:metrizability}
  By Proposition~\ref{prop:kR-space}, if $\cP$ is countable, then the topology on~$L$ is metrizable.
  Even when $\cP$ is countable, however, our purposes are better suited by thinking of~$L$ as endowed with the topology (and corresponding uniformity~\cite[§~8.1]{Engelking1989}) explicitly given by the full predicate collection, rather than by an implied “master” metric which, in an abridged manner, induces the same topology.
\end{remark}

\subsection{Tychonoff and Realcompact spaces}
\label{sec:Tychonoff-RC}

\subsubsection{Tychonoff spaces}
\label{sec:Tychonoff}
Recall that a topological space~$X$ is \emph{Tychonoff} if it is $T_{3\text{½}}$, i.e., a completely regular Hausdorff space;
explicitly: (\emph{i})~points are closed, and (\emph{ii})~given any point~$x\in X$ and closed~$C\not\ni x$ there exists a continuous function~$f: X\to[0,1]$ such that $f(x)=0$ and $f{\restriction}C=1$.

\begin{remark}\label{rem:CSS-def-undef-preds}
  A reduced CSS $\underline{L} = \langle L, \cP\rangle$ is ultimately just a Tychonoff space endowed with a distinguished family~$\cP$ of real functions $P(\cdot)$ (distinguished predicates), and such that the topology on~$L$ is initial by the collection~$\cP$
  (i.e., the topology of~$L$ is generated by the inverse images of open intervals of~$\RR$ under observables~$P(\cdot)$).
  From another perspective, any CSS $\underline{L} = \langle L, \cP\rangle$ is isomorphic to a “sub-CSS” of a CSS $\underline{\RR^{\cP}} = \langle\RR^{\cP}, (\pi_P)_{P\in\cP}\rangle$ (where $\pi_{P}$ is the projection $v\mapsto v_P$) via the type map $\tp: L\to\RR^{\cP}$, which is a homeomorphic feature-preserving embedding;
  therefore, such product CSSs $\underline{\RR^{\cP}}$ are universal.

  The observable features~$P(\cdot)$ are regarded as being “computable on~$L$” \emph{ab initio};
  they also may be seen as monomials generating a polynomial algebra of continuous real functions on~$L$;
  the closure $\mathbf{D}\subseteq\RR^L$ of the set of such polynomials in the topology of uniform convergence on compacta $K\subseteq L$ (equivalently, on shards~$L[\rb]$) is called the algebra of \emph{definable predicates} on~$L$.
  If $L$ is compact, $\mathbf{D}=\mathrm{C}(L)$ is the set of all continuous functions by the Stone--Weierstrass Theorem;
  if $\cP$ is countable, then $\mathbf{D}\subseteq \mathrm{C}(L)$
  (cf., Proposition~\ref{prop:kR-space}(4));
  in general, however, $\mathbf{D}$ is neither contained in nor a superset of~$\mathrm{C}(L)$.
  Any non-definable predicate~$\varphi$ is “transcendental” over~$\cP$
 —not merely in an algebraic sense, but in a stronger topological one:
  not only does such a non-definable $\varphi$ fail to be a polynomial on monomials~$P(\cdot)\in\cP$:
  it is not even possible to approximate $\varphi$ by such polynomials, uniformly on shards.
\end{remark}

\subsubsection{Realcompact spaces}
\label{sec:realcompact}
A topological space is called \emph{realcompact} if it is Tychonoff and it embeds homeomorphically as a \emph{closed} subspace of the topological product $\RR^I = \prod_{i\in I}\RR$ for some index set~$I$~\cite[§~3.11]{Engelking1989}.
(There is a multitude of equivalent definitions of realcompactness.
For a thorough treatment of realcompact spaces, refer to Weir’s monograph~\cite{Weir75}.)

A CSS $\underline{L} = \langle L,\cP\rangle$ is realcompact iff the type map $L\to\RR^{\cP}: v\mapsto\tp(v)\coloneqq(P(v))_{P\in\cP}$ has \emph{closed} image $\tp(L)\coloneqq\{\tp(v): v\in L\}$ in~$\RR^{\cP}$, i.e., if $\tp(L) = \fL$ is the full space of state types of~$L$
(all state types are realized).
Any compact (Hausdorff) CSS~$\underline{L}$ is necessarily Tychonoff and in fact realcompact:
Taking $\cP$ to be any set of continuous functions $P: X\to\RR$ separating points of~$X$, the type map $\tp: X\to\RR^{\cP}$ is injective and has compact, hence closed, image;
it is therefore a homeomorphic embedding.

\subsubsection{Realcompactness of type spaces}
\label{sec:fL-is-RC}

Any type space $\fL \subseteq \RR^{\cP}$ is a closed subspace, hence realcompact.
Identifying the layer space $L$ with its embedded image $\tp(L)\subseteq\fL$, it is suggestive to regard the realcompact type space~$\fL$ as a canonical realcompact extension of~$\tp(L)\cong L$.
Such viewpoint is quite appropriate for our purposes, so we discuss in what precise sense this realcompact extension is canonical.

More generally, consider any Tychonoff space~$X$ whose topology is initial with respect to a collection $\Phi$ of real functions $\varphi: X\to\RR$
(i.e., inverse images of opens of~$\RR$ by such functions generate the topology of~$X$).
Each point $x\in X$ has a \emph{$\Phi$-type} $\Ftp(x) = (\varphi(x))_{\varphi\in\Phi}\in\RR^{\Phi}$, and $X$ embeds (via the map $\Ftp$) as a subspace $\Ftp(X)\subseteq\RR^{\Phi}$ whose closure $\fX = \overline{\Ftp(X)}\subseteq\RR^{\Phi}$ (the \emph{$\Phi$-type space of~$X$}) is realcompact.

The type space $\fX = \fX_{\Phi}$ depends on~$\Phi$.
A key observation is that each of the functions $\varphi\in\Phi$ extends to~$\fX$ continuously (as the “$\varphi$-th coordinate” map $\hat{\varphi}: (\fx_{\psi})_{\psi\in\Phi}\mapsto \fx_{\varphi}$).
However, other real functions on~$X$ —even if continuous— need not extend to~$\fX$ continuously.
Thus, the $\Phi$-type space $\fX \supseteq \Ftp(X)\cong X$ possesses the universal property that every $\varphi\in\Phi$ admits a unique continuous extension to~$\fX$;
it is characterized by such universal property up to homeomorphism.%
\footnote{In fact, any real function~$\xi$ on~$X$ that is uniformly approximable by polynomials in the functions $\varphi\in\Phi$ is (necessarily continuous, and) extends continuously to a real function $\hat{\xi}$ on~$\fX$
(a uniform limit of polynomials in the corresponding functions~$\hat{\varphi}$ on~$\fX$),
so $\fX$ possesses the extension property for functions in the uniform closure of the real algebra generated by functions~$\varphi\in\Phi$.}

\begin{remark}\label{rem:RC-CCS}
Let $C\coloneqq \mathrm{C}(X)$ be the set of \emph{all} continuous real-valued functions on a Tychonoff space~$X$, and let $\hat{X}\subseteq\RR^{C}$ be the corresponding type space, called the \emph{realcompactification of~$X$}.
  Every continuous function $\varphi$ on~$X$ extends to a continuous function~$\hat{\varphi}$ on the realcompactification~$\hat{X}$.
  In fact, for any $\Phi\subseteq C = \mathrm{C}(X)$, the $\Phi$-type space~$\fX\subseteq\RR^{\Phi}$ is a quotient (not a subspace!)\ of $\hat{X}\subseteq\RR^{C}$ in a natural manner:
  indeed, $\fX\subseteq\RR^{\Phi}$ is the image of~$\hat{X}\subseteq\RR^{C}$ under the natural projection map $\RR^C\to\RR^{\Phi}: (\fx_{\varphi}: \varphi\in C)\mapsto (\fx_{\varphi}: \varphi\in\Phi)$.
  Thus, given a fixed set~$\Phi$ of continuous real functions on~$X$, it is appropriate to think of the $\Phi$-type space~$\fX$ as a “$\Phi$-relative realcompactification” of~$X$, since $\fX\supseteq X$ possesses the universal extension property \emph{only} for functions~$\varphi\in\Phi$
 —rather than for \emph{all} $\varphi\in\mathrm{C}(X)$, which corresponds to the (“absolute”) realcompactification~$\hat{X}$ of~$X$.%
 \footnote{The notation~$\upsilon X$ (“upsilon-$X$”) is standard for the realcompactification —denoted $\hat{X}$ above— of a Tychonoff space~$X$.}

\end{remark}

\subsubsection{Realcompact CSSs}
\label{sec:RC-CSSs}
We single out the (sub)class of CSSs possessing the:
\begin{itemize}
\item \emph{\textbf{Realcompactness Property.}\quad Every type $\fv\in\fL$ is realized, i.e., of the form $\tp(v)$ for some $v\in L$.}
\end{itemize}
Thus, realcompactness is the requirement that $\tp: L\to \fL$ be a \emph{surjection} onto the type space~$\fL$, whence $\tp$ is a homeomorphism $L\cong\fL$ (by the Reduction Axiom).
Rephrasing, realcompactness states that whenever $(u^{(i)})_{i\in I}\subseteq L$ and $\cU$ on $I$ are such that the ultralimit $\fv\coloneqq\Ulim_i\tp(u^{(i)})$ exists, then some $v\in L$ satisfies $\tp(v) = \fv$.

It is appropriate to regard realcompactness as capturing a certain notion of “completeness” or “saturation” of the space~$L$.
Particularly when $\cP$ is infinite, realcompactness is a rather strong requirement on CSSs, so we do not impose it as an axiom;
instead, we rely primarily on the realcompactness and universal properties of the type space~$\fL\supseteq L$.%
\footnote{When $\cP$ is finite, realcompactness is a rather mild requirement:
  it is seen to be equivalent to the completeness of $\langle L,\delta\rangle$, where $\delta$ is the metric introduced in the proof of Proposition~\ref{prop:kR-space}.}

\subsection{Compositional Computation Structures}
\label{sec:CCS}

A \emph{Compositional Computation Structure (CCS)}
\begin{equation*}
  \cC = \langle \underline{L},\underline{\Gamma},\ev\rangle
\end{equation*}
with observables collection $\cP$ consists of
\begin{itemize}
\item a CSS $\underline{L} = \langle L, (P(\cdot))_{P\in\cP}\rangle$ with predicate symbol set~$\cP$ and, for each $P\in\cP$, an observable $P(\cdot):L\to\RR$;
\item a semigroup $\underline{\Gamma} = (\Gamma,\circ)$, the \emph{computations sort}
  (the —associative— semigroup operation~$\circ: \Gamma\times\Gamma\to\Gamma$ is denoted simply $(\gamma,\delta)\mapsto \gamma\delta$ when convenient);
\item a map $\ev: \Gamma\times L\to L: (\gamma,v)\mapsto \ev(\gamma,v)$ (\emph{“evaluation”}) giving an action of $\Gamma$ on~$L$
  (i.e., $\ev(\gamma\delta,v) = \ev(\gamma,\ev(\delta,v))$ for $\gamma,\delta\in\Gamma$ and $v\in L$).
\end{itemize}
\begin{remarks}\label{rem:n-ary-transfns}
  \begin{enumerate}
  \item In principle, the semigroup operation ‘$\circ$’ of~$\underline{\Gamma}$ and evaluation action ‘$\ev$’ are abstract
    (i.e., not literally composition and application of functions).
    However, one may always regard $\underline{\Gamma}$ as a semigroup (under the operation ‘$\circ$’ interpreted as composition) of maps $\gamma(\cdot): v\mapsto\ev(\gamma,v)$;
   —i.e., regard $\Gamma$ as a sub-semigroup of the semigroup $L^L$ of all maps $L\to L$, under \emph{bona fide} functional composition:
 nothing of structural relevance is lost thus.
 The structural viewpoint abstracts inessential aspects of a concrete such realization of~$\Gamma$.
 In practice, it is convenient to identify $\gamma\in\Gamma$ with $\gamma(\cdot)\in L^L$.
  \item In applications, a more general notion of \emph{CCS with $n$-ary computations} is useful.
    By this we mean that computations $\gamma\in\Gamma$ may each have an \emph{arity} $n = n_{\gamma} \in \NN$ such that $\gamma(\cdot)$ is an ($n$-argument) map $L^n\to L$.
    (It is appropriate to regard evaluation on $n$-ary such $\gamma$ as a map $\ev_n:\Gamma_{\!n}\times L^n\to L$, where $\Gamma_{\!n}\subseteq\Gamma$ is the set of $n$-ary elements~$\gamma$;
    thus,  $\gamma\in\Gamma_{\!n}$ gives an $n$-argument map $\ev(\gamma;\cdot): L^n\to L$.)
    CCS with $n$-ary computations augments the semigroup operation ${\circ}: \Gamma_1\times\Gamma_1\to\Gamma_1$ with a richer set of operations realizing arity-appropriate compositions.
    Explicitly, given $\gamma_1,\gamma_2,\dots,\gamma_m\in \Gamma_{\!n}$ and $\theta\in\Gamma_{\!m}$, there exists an element $\eta\coloneqq \theta{\circ}(\gamma_1,\dots,\gamma_m)\in \Gamma_{\!n}$ satisfying
    \begin{equation*}
      \ev(\eta;\bar{v}) =
      \ev\bigl(\theta;
      \ev(\gamma_1,\bar{v}),
      \ev(\gamma_2,\bar{v}),
      \dots,
      \ev(\gamma_m,\bar{v})
      \bigr)
      \qquad
      \text{for all $\bar{v} \in L^n$,}
    \end{equation*}
    i.e., the above identity holds for a suitable “generalized composition” operation~$\circ$
   —or, rather, for one such operation ${\circ}_m^n: \Gamma_{\!m}\times(\Gamma_{\!n})^m\to\Gamma_{\!n}$ for each $m,n\in\NN$—
    moreover, the sort of computations $\underline{\Gamma} = (\Gamma,\circ_m^n: m,n\in\NN)$ is endowed with all such compositions.
  \end{enumerate}
    In order to simplify the exposition, CCSs with $n$-ary computations as in the preceding remark will be used only in informal discussions and examples.
\end{remarks}

\subsubsection{Reduction and Continuity Axioms}
\label{sec:cont-redn-axioms}

Every CCS $\cC = \langle\underline{L}, \underline{\Gamma}\rangle$ will be assumed to satisfy the following axioms:

\begin{itemize}
\item \emph{\textbf{Reduction Axioms for Compositional Computation Structures.}
    \begin{enumerate}
    \item States $v,w\in L$ are equal if their types $\tp(v), \tp(w)$ are equal
      (i.e., the underlying CSS~$\underline{L}$ is reduced);
    \item Transformations $\gamma,\delta\in\Gamma$ are equal if the maps $\gamma(\cdot), \delta(\cdot): L\to L$ are equal.
    \end{enumerate}
  }
\end{itemize}

As a temporary (weaker) placeholder for the Extendability Axiom (see~§~\ref{sec:Ax-Extend}) eventually imposed on CCSs, we presently impose the natural:
\begin{itemize}
\item \emph{\textbf{Continuity Axiom:}\quad
    The action of $\underline{\Gamma}$ on~$L$ is by maps continuous in the topology of~$L$
  (i.e., is a topological action on the CSS $\underline{L}$)}.
\end{itemize}
Explicitly, for each computation $\gamma\in\Gamma$ and~$P\in\cP$, the real-valued “$P$-feature” $P{\circ}\gamma: v\mapsto P(\gamma(v))$ of~$\gamma(\cdot)$ is continuous on~$L$.

Reduction Axiom~(2) says that $\Gamma$ is bijectively identified with its image $\Gamma(\cdot)\subseteq L^L$
of maps (“state transitions”) $\gamma(\cdot)$ ($\gamma\in\Gamma$).
This identification gives a natural topology on~$\Gamma$, obtained (as pullback) from the topology of pointwise convergence on the maps $\gamma(\cdot)\in\Gamma(\cdot)\subseteq L^L$ associated to computations~$\gamma\in\Gamma$; the Reduction Axiom implies that this topology is also Hausdorff.

As long as the Continuity Axiom holds, the Reduction Axioms are innocuous requirements on a CCS~$\cC$, because one can always pass from $\cC$ to a reduced CCS~$\tilde{\cC}$ (i.e., one satisfying the Reduction Axiom) as follows.
First, replace $\Gamma$ by its quotient $\tilde{\Gamma}\coloneqq \Gamma/(\tp{\circ}\ev)$ upon identifying computations $\gamma,\delta\in\Gamma$ such that $\tp(\ev(\gamma,v)) = \tp(\ev(\delta,v))$ for all $v\in L$;
second, pass from the underlying CSS $\underline{L}$ to its quotient-by-type $\underline{\tilde{L}}$ if necessary.
The evaluation $\ev$ of~$\cC$ induces a well-defined natural action $\widetilde{\ev}: \tilde{\Gamma}\times \tilde{L} \to \tilde{L}$.
By the Continuity Axiom, the passage from $\cC$ to $\tilde{\cC} = \langle \underline{\tilde{L}}, \underline{\tilde{\Gamma}}, \widetilde{\ev}\rangle$ preserves all structural properties of states and computations, as well as the Continuity Axiom.%
\footnote{The passage to $\tilde{\cC}$ also preserves the Extendability Axiom of §~\ref{sec:Ax-Extend}.}

\begin{remark}\label{rem:non-conts-extn}
  The Continuity Axiom ensures that computations act continuously on~$L$.
  In general, however, the action $\gamma(\cdot)$ of a computation $\gamma\in\Gamma$ on the state space~$L$ need \emph{not} admit a continuous extension to a transformation $\fL\to\fL$.
  This distinction is quite important;
  it speaks to the weakness of the Continuity Axiom, and suggests a strengthening
  —namely, the Extendability Axiom below, which is a key assumption in our main results. 
\end{remark}

\subsection{Examples of CSSs and CCSs}
  \label{sec:CCS-exa}

\subsubsection{The unit interval}
\label{sec:unit-interval}

Consider a CSS with state space $L = [0,1]$ (the unit interval) endowed with the single identity predicate $\id: [0,1]\to\RR:v\mapsto v$.
This CSS $\langle[0,1],\id\rangle$ is realcompact.%
\footnote{We write simply “$\id$” for the observables collection of this CSS rather than the cumbersome “$(\id)$” denoting a singleton collection of observables.}

Let $\underline{\Gamma} = (\Gamma,\circ)$ be any semigroup (under composition) of continuous functions $\gamma: [0,1]\to [0,1]$, acting on~$[0,1]$ by functional application: $\ev(\gamma,v) \coloneqq \gamma(v)$;
this yields a realcompact CCS $\langle([0,1],\id),\underline{\Gamma}\rangle$.
An interesting such CCS has semigroup $\Gamma = \{\gamma^n: n\in\NN\}$ consisting of iterates of the chaotic map $\gamma: v\mapsto 4v(1-v)$.

Replacing $[0,1]$ with the open interval $(0,1)$, one obtains a non-realcompact CSS $\langle(0,1), \id\rangle$ having (realcompact) type space $\overline{(0,1)} = [0,1]\subseteq\RR^1 = \RR$.
 (By contrast —cf., Remark~\ref{rem:RC-CCS}— the realcompactification $\widehat{(0,1)} \supseteq (0,1)$ is a much larger topological extension \emph{not} homeomorphic to a subset of~$\RR$.)

\subsubsection{$\RR^d$}
\label{sec:RR_d}

Given $d\in\NN$, we obtain a CSS $\underline{\RR^d} = \langle\RR^d,(P_i)_{i=1}^d\rangle$ on the $d$-dimensional real space $\RR^d$ (of states), endowed with coordinate functions $P_i(v) \coloneqq v_i$ ($1\le i\le d$) as observable features.
The corresponding type space is $\tp(\RR^d) = \RR^d$;
the type topology coincides with the usual one, so $\underline{\RR^d}$ is realcompact.

    There is ample flexibility in expanding the collection of observable features with new predicates yielding formally distinct CSSs $\underline{L}$ with layer states sort~$L = \RR^d$.
  For any real $q\ge 1$, one may (for instance) expand the predicate collection $\cP$ with the “$q$-norm” predicate $\Nrm{\cdot}_q$ defined by
  \begin{equation*}
    \Nrm{v}_q \coloneqq \sqrt[q]{\nrm{v_1}^q + \dots + \nrm{v_d}^q}.
  \end{equation*}
  One may also expand the observables collection with, say, the supremum norm
  \begin{equation*}
    \Nrm{v}_{\sup} = \Nrm{v}_{\infty} \coloneqq \sup_{1\le i\le d}\nrm{v_i}\quad
    \bigl(= \max(\nrm{v_1}, \dots, \nrm{v_d})\bigr).
  \end{equation*}
  Since $d$ is finite, the predicates $\Nrm{\cdot}_q$ above are continuous with respect to the topology of~$\RR^d$.
  In fact, any continuous predicate $\varphi:\RR^d\to\RR$ may be added to the collection of observables of $\underline{\RR^{d}}$ yielding an essentially equivalent CSS:
  any such CSS is an expansion of $\underline{\RR^{d}}$ by \emph{definable predicates} in the sense of §~\ref{sec:definable-predicates} below;
  moreover, it is homeomorphic to $\underline{\RR^{d}}$ (still realcompact).%
  \footnote{On the other hand, whether an arbitrary CSS is realcompact or not, any expansion thereof by a \emph{bounded} predicate that is \emph{dis}continuous results in a CSS that is not realcompact.
    (Boundedness —or at least boundedness near points of discontinuity— cannot be dispensed with: expanding $\underline{\RR}$ by the unbounded, discontinuous predicate $\varphi$ given by $\varphi(t)\coloneqq 1/t$ for $t\ne0$ and $\varphi(0)\coloneqq 0$ yields a \emph{realcompact} CSS, since the resulting set of realized types —the graph of~$\varphi$— is closed in~$\RR^{2}$.)}

Expanding $\underline{\RR^d}$ with:
  \begin{itemize}
    \item $\overline{\Gamma} = (\Gamma,\circ)$ any semigroup of continuous functions on~$\RR^d$; and
    \item the evaluation action of~$\Gamma$ on~$L$ by functional application $\ev(\gamma,v)\coloneqq \gamma(v)$ as in~\ref{sec:unit-interval} above,
    \end{itemize}
    one obtains a CCS $\langle\RR^d,(P_i)_{i=1}^d,\overline{\Gamma}\rangle$.
    Natural such expansions are by semigroups $\Gamma$ of linear or affine operators on~$\RR^d$.

\subsubsection{$\RR^{\omega}$ and $c_{00}$}
\label{sec:RR_omega}

Let the CSS $\underline{\RR^{\omega}}$ have states space $\RR^{\omega}\coloneqq \prod_{i\in\NN}\RR$ consisting of all real sequences~$v = (v_i)_{i\in\NN} \subseteq \RR$, endowed with the collection~$\cP$ of predicates $P_i: v\mapsto v_i$ ($i\in\NN$).
Such CSS $\underline{\RR^{\omega}}$ is realcompact.

The subspace $c_{00}$ consisting of sequences $v$ having at most finitely many entries $v_i\ne 0$ is a non-realcompact sub-CSS of~$\underline{\RR^{\omega}}$ whose type space is homeomorphic to the product space~$\RR^{\omega}$.

A natural expansion of $\RR^{\omega}$ to a CCS is by the semigroup $\underline{\Gamma}$ of continuous linear operators thereon.
Each linear such computation $\gamma\in\Gamma$ is effectively a collection $(\gamma_i)_{i\in\NN}$ of real functionals $\gamma_i \coloneqq P_i{\circ}\gamma: \RR^{\omega}\to\RR$, each of the form $\gamma_i: v\mapsto \sum_{j\in\NN}r^{(i)}_jv_j$, for some scalar collection $\rb^{(i)} = (r^{(i)}_j)_{j\in\NN}\in c_{00}$.
Thus, every entry $w_i = P_i(\gamma(v))$ of $w = \gamma(v)$ is a \emph{finite} linear combination of features $P_j(v)$ of the input~$v$.%
\footnote{The topological vector spaces $c_{00}$ and $\RR^{\omega}$ are topological duals of one another, so one obtains a corresponding (dual) characterization of continuous linear maps on~$c_{00}$.}

Many natural real functions on~$c_{00}$ (or on suitable subspaces thereof) are discontinuous (in the topology of entry-wise convergence);
expanding $c_{00}$ with any such function as distinguished predicate leads to (non-homeomorphic) CSSs.
  
 \begin{remark}
 Note that in the CSSs~\ref{sec:unit-interval}–\ref{sec:RR_omega} above (but not in~\ref{sec:ell_q} below), a state~$v$ is exactly the same as its type~$\tp(v)$.
 \end{remark}

 \subsubsection{$\ell_q$}
\label{sec:ell_q}

 For any extended real $q\in[1,\infty]$, consider the states space
  \begin{equation*}
    L = \ell_{q} = \{v\in\RR^{\omega}: \Nrm{v}_{q}<\infty\}.%
    \footnote{The norm $\Nrm{\cdot}_q$ is the classical one as in~\ref{sec:RR_d} above, namely $\Nrm{v}_q \coloneqq \bigl(\sum_i\nrm{v_i}^q\bigr)^{1/q}$ for $q<\infty$, and $\Nrm{v}_{\infty} \coloneqq \sup_i\nrm{v_i}$.}
  \end{equation*}
  For $q<\infty$, such space $\ell_q$ is the $\Nrm{\cdot}_q$-metric completion of the subspace $c_{00}\subseteq\RR^{\omega}$ and thus norm-complete (i.e., a Banach space);
  $\ell_{\infty}$ is also norm-complete.%
  \footnote{The $\Ninf{\cdot}$-metric completion of $c_{00}$ is the separable space $c_0 = \{v\in\RR^{\omega}: \lim_{i\to\infty}v_i=0\} \subsetneq \ell_{\infty}$
   —the space~$\ell_{\infty}$ is not separable.}
  A natural collection of observables is $\cP = (N_q, P_i)_{i\in\NN}$, where $P_i:v\mapsto v_i$ is the $i$-th coordinate as in~\ref{sec:RR_omega} above, and $N_q$ names the norm $N_q(\cdot): v\mapsto\Nrm{v}_q$.
  Since the predicate collection $\cP$ is countable, the type topology on the layer space $\ell_q$ is metrizable;
  for $q<\infty$ it coincides with the usual $\Nrm{\cdot}_q$-norm topology, whereas for $q=\infty$ it is strictly coarser.%
  \footnote{Cf., the proof of Proposition~\ref{prop:standard-types-countable-P} below for metrizability.
    For $q<\infty$: given $v^{(n)}\to v$ entry-wise with $\Nrm{v^{(n)}}_q\to\Nrm{v}_q$, one bounds $\Nrm{v^{(n)}-v}_q$ upon splitting the sum at an index~$N$ beyond which the tail $\sum_{i>N}\nrm{v_i}^q$ is small (the hypotheses force the corresponding tails of the~$v^{(n)}$ to be small as well).
    For $q=\infty$: the sequence $v^{(n)}\coloneqq e_1+e_n$ ($n\ge2$) —where $e_i$ denotes the $i$-th unit vector— converges to~$e_1$ in the type topology (entry-wise, and $\Ninf{v^{(n)}} = 1 = \Ninf{e_1}$), yet $\Ninf{v^{(n)}-e_1} = 1$ for all~$n$.}
  In either case, $\ell_q$ is \emph{non}-realcompact.
  It is easy to see that its type space is $\fL_q = \{(v,r)\in\ell_q\times\RR: r \ge \Nrm{v}_q\}$.
    (The set of realized types is $\tp(\ell_q) = \{(v,r)\in\ell_q\times\RR: r = \Nrm{v}_q\} \subseteq \fL_q$, for which the “correct” norm $\Nrm{v}_q$ agrees with the interpretation value $N_q(v,r) = r$ of the symbol~$N_q$.)
    Fixing $q<\infty$, the function $\Ninf{\cdot}: \ell_q\to\RR: v\mapsto \sup_n\nrm{P_n(v)}$ is $1$-Lipschitz, hence continuous on~$\ell_q$;
    however, the corresponding function $\tp: \ell_q\to\RR: \tp(v)\mapsto\Ninf{v}$ does not extend continuously to~$\fL_q$.%
    \footnote{By the Stone--Weierstrass Theorem, every continuous function on the compact Hausdorff space $\fL_q[1] \coloneqq \{(v,r)\in\fL_q: r\le 1\}$ is uniformly approximable by algebraic combinations (finitely many at a time) of predicates $P_i$, and~$\Nrm{\cdot}_q$;
      however, an elementary argument shows that $\tp(v)\mapsto\Ninf{v}$ admits no such uniform approximations on~$\fL_q[1]$.}

    A natural expansion of $\underline{\ell_q}$ to a CCS is by its semigroup $\underline{\Gamma}$ of bounded (i.e., norm-continuous) linear operators.
    Such operators are continuous on~$\underline{\ell_q}$;
    however, they are typically discontinuous when regarded as functions on the reduct CSS of~$\underline{\ell_q}$ wherein the additional predicate $N_q$ is removed, i.e., when $\ell_q$ is topologized as sub-CSS of~$\RR^{\omega}$ rather than of $\RR^{\omega}\times\RR$ as above.

\begin{remark}
  The metric $\dd_q\coloneqq (v,w) \mapsto\Nrm{w-v}_q$ on $\ell_q$ (or on $\cL^d$ for $d$ finite) is not definable in terms of the norm predicate $\Nrm{\cdot}_q$ unless $\ell_q$ is expanded to a CCS with, say, the \emph{binary} predicate of subtraction $(v,w)\mapsto v-w$.
  While otherwise not meant to detract from the preceding discussion, this remark serves to highlight the usefulness of CCSs with $n$-ary layer transformations (cf., Remark~\ref{rem:n-ary-transfns}).
\end{remark}

\section{Deep computations}
\label{S:deep-comp}

Throughout this section, $\cC = \langle \underline{L},\underline{\Gamma}\rangle$ will be a fixed CCS.

\subsection{Shards in state- and type-spaces}\label{sec:cpct-type-spaces}

\subsubsection{Sizers and shards in type spaces}
\label{sec:sizers}
A \emph{sizer} is any collection $\rb = (r_P)_{P\in\cP}\in[0,\infty)^{\cP}$ of nonnegative reals.
The number $r_P$ is called an \emph{a priori} bound for~$P$.

For a sizer~$\rb$, we introduce the topological product space
\begin{equation*}
\RR{[\rb]} \coloneqq \prod_{P\in\cP}[-r_P,r_P];
\end{equation*}
it is a compact subspace of the product space~$\RR^{\cP}$;
moreover, $\RR^{\cP} = \bigcup_{\rb}\RR{[\rb]}$ (with $\rb$ varying over sizers).
A subset $S\subseteq\RR^{\cP}$ is called \emph{entry-wise bounded} if $S\subseteq\RR{[\rb]}$ for some sizer~$\rb$.
Clearly, \emph{relatively compact subsets of $\RR^{\cP}$ are precisely entry-wise bounded subsets}.

\subsubsection{Shards}
\label{sec:shards}

For a sizer~$\rb$, the \emph{$\rb$-shard of~$L$} is
\begin{equation*}
  L[\rb] \coloneqq \{v\in L: \tp(v)\in\fL[\rb]\}
  = \{v\in L: \text{$\nrm{P(v)}\le r_P$ for all $P\in\cP$}\}.
\end{equation*}
Clearly, an arbitrary intersection of shards is a shard, and any finite union of shards is included in some shard.

Let the \emph{type $\rb$-shard} $\fL[\rb]\subseteq\RR{[\rb]}$ be the topological closure of the set $\tp(L[\rb]) \coloneqq \{\tp(v): v\in L[\rb]\}$ of types realized in (i.e., by elements of)~$L[\rb]$.
Evidently, $\fL[\rb]\subseteq \fL\cap\RR{[\rb]}$;
however, the inclusion is typically strict because a type $\fu\in\fL\cap\RR{[\rb]}$ need not be an accumulation point of types realized in $L[\rb]$ itself, thus $\fu$ need not belong to~$\fL[\rb]$.%
\footnote{In general, a type $\fu\in\fL\cap\RR{[\rb]}$ need not even be an accumulation point of realized types in \emph{any} fixed shard~$L[\sbb]$!}
We introduce the space $\LSh \coloneqq \bigcup_{\rb}\fL[\rb]$ of \emph{shard-supported types};
it is the set of types~$\fv$ of arbitrary shards~$L[\rb]$.
By the preceding discussion, we have $\LSh \subseteq\fL$, but the inclusion is proper in general
($\LSh $ need not be closed in~$\fL$).
The space $\LSh $ will be of central importance in what follows.

A collection~$R$ of sizers~$\rb$ is \emph{exhaustive} if, for any sizer $\sbb$ there exists $\rb\in R$ such that $r_P\ge s_P$ for all $P\in\cP$.

From its definition, it is clear that $\LSh $ is the union of type-shards $\fL[\rb]$ as $\rb$ varies over any exhaustive~$R$.

Recall that a Hausdorff space $X$ is
\begin{itemize}
\item a \emph{k-space} if closed subsets $Y\subseteq X$ are precisely those whose intersection $Y\cap K$ with an arbitrary compact subset $K\subseteq X$ is closed~\cite[3.3.18ff]{Engelking1989};
\item a \emph{\kR-space} if an arbitrary real function $\varphi:X\to\RR$ is continuous as soon as its restrictions $\varphi{\restriction}K$ to compacta $K\subseteq X$ are continuous.
\end{itemize}
Evidently, any k-space is a \kR-space.

\begin{proposition}\label{prop:standard-types-countable-P}
  Let $\underline{L} = \langle L,\cP\rangle$ be a CSS whose distinguished predicate collection $\cP = (P_i)_{i\in\NN}$ is countable.
  \begin{enumerate}
  \item $\LSh  = \fL$ (i.e., all types are shard-supported).  
  \item $\LSh$ is a k-space.
    More precisely, closed subspaces $S\subseteq\LSh$ are (precisely) those whose intersection $S\cap\fL[\rb]$ with an arbitrary type-shard $\fL[\rb]$ is closed.
  \end{enumerate}
\end{proposition}
\emph{A fortiori}, the result holds when the predicate collection is finite.

 We thank F.\ Tall for bringing to our attention the fact that realcompact spaces embeddable in $\RR^{\omega}$ are k-spaces.
\begin{proof}
  For $i<\omega$, introduce the pseudometric $\dd_i(\fu,\fv) \coloneqq \nrm{P_i(\fu)-P_i(\fv)}$ on~$\fL$, and let $\delta_i\coloneqq \dd_i/(1+\dd_i)$ be the usual $[0,1]$-valued pseudometric corresponding to~$\dd_i$.
  The space~$\fL$ is completely metrizable by $\dd(\fu,\fv) \coloneqq \sum_{i\in\NN}2^{-i}\delta_i(\fu,\fv)$, and the subset $\tp(L)\subseteq\fL$ of realized types is dense.
  \begin{enumerate}
  \item 
  Given $\fu\in\fL$, there is a sequence $(v_n)_{n\in\NN}$ such that $\lim_{n\to\infty}\tp(v_n)=\fu$ in the $\dd$-metric sense.
  The set $A\coloneqq\{\tp(v_n): n\in\NN\}\cup\{\fu\}\subseteq\fL$ is compact
  (any open cover $\mathcal{O}$ of~$A$ contains an open $U\ni \fu$;
  since $\tp(v_n)\to \fu$, all but finitely many elements of~$A$ are contained in~$U$, so $\mathcal{O}$ has a finite subcover).
  By compactness of~$A$, the image $P(A)\coloneqq \{P(\fv): \fv\in A\}$ is compact in~$\RR$ for each $P\in\cP$, hence bounded, say $\nrm{P(\fv)}\le r_P$ for some $r_P\ge0$ and all $\fv\in A$.
  Clearly, $\fL[\rb]\supseteq A\ni \fu$, so $\fL\subseteq\LSh $ ($\subseteq\fL$, in any case);
  hence, $\LSh  = \fL$.

\item %
  If $S \subseteq \fL$ ($=\LSh$) is not closed, let $\fu\in \overline{S}\setminus S\subseteq \fL\setminus S$.
  As in~(1) above, construct $(\fv_n)_{n\in\NN}\subseteq S$ with $\fu = \lim_{n\to\infty}(\fv_n)$.
  Since $\overline{L} = \fL \supseteq S$, for each $n,k\in\NN$ there exists some $v^{(n)}_k\in L$ such that $\dd(v^{(n)}_k, \fv_n)<1/(n+k)$.
  One sees that $A = \{\tp(v^{(n)}_k):n,k\in\NN\} \cup \{\fv_n:n\in\NN\} \cup \{\fu\}$ is compact.%
  \footnote{Clearly, $A^{(n)}\coloneqq \{\fv_n\} \cup \{\tp(v^{(n)}_k): k\in\NN\}$ is compact for each $n\in\NN$ (as in~(1) above).
    Given any open cover $\mathcal{O}$ of~$A$, some $U\in \mathcal{O}$ satisfies $U\ni \fu$.
    For a set $N\subseteq\NN$ containing all but finitely many~$n$, we have $U\supseteq A^{(n)}$.
    The set $B\coloneqq \bigcup_{n\in\NN\setminus N} A^{(n)}$ is a finite union of compacta, hence itself compact.
    Thus, finitely many opens of~$\mathcal{O}$ cover~$B$, which together with~$U$ cover~$A$.}
  Therefore, $A\subseteq\fL[\rb]$ for some sizer~$\rb$, and thus $\{v^{(n)}_k:n,k\in\NN\}\subseteq L[\rb]$, hence $\fu = \lim_{n\to\infty}\tp(v^{(n)}_n)\in \overline{\tp(L[\rb])} = \fL[\rb]$.
    Since $\{\fv_n: n\in\NN\}\subseteq S\cap\fL[\rb]$ and its limit~$\fu$ lies in $\overline{S}\setminus S$, we see that $S\cap\fL[\rb]$ is not closed. %
    The converse is trivial, since type-shards $\fL[\rb]$ are closed
    (and the intersection of closed sets is closed).%
    \footnote{We recall the following fact, closely related to~(2):
      Every sequential Hausdorff (therefore, every metric) space is a k-space~\cite[Theorem~3.3.20]{Engelking1989}.}\qedhere
\end{enumerate}
\end{proof}

\begin{remarks}\label{rem:cofinal-sizer-set}
  \begin{enumerate}
  \item For $\cP$ at most countable, Proposition~\ref{prop:standard-types-countable-P} implies that
  closed subsets $S\subseteq\LSh$ are those whose intersections $S\cap\fL[\rb]$ are closed for all sizers~$\rb$ in any exhaustive collection~$R$.
\item A compact subset $K\subseteq\LSh$ need not be included as a subset of any type-shard~$\fL[\rb]$ (even if $\cP$ is countable).
  \end{enumerate}
\end{remarks}

\subsection{Transitions-in-type. Extendability.}\label{sec:t-t}

Any —not necessarily continuous— function $f:L\to L$ (resp.,  $f: L\to\LSh$) %
will be called a \emph{(layer) transition}
(resp., an \emph{ultra-transition (u-t)}).
We also introduce the notion of \emph{transition-in-type (t-t)} to mean any function $\ff: \LSh\to\LSh$.

A transition $f$ (resp., a u-t~$f$; a t-t~$\ff$) is \emph{shard-to-shard (Sh2Sh)} if, for every sizer~$\rb$ there is a sizer $\sbb$ such that $f$ restricts to a map~$L[\rb]\to L[\sbb]$
(resp., $f$ restricts to~$L[\rb]\to\fL[\sbb]$;
$\ff$ restricts to $\fL[\rb]\to\fL[\sbb]$).

The \emph{transition space}, \emph{ultra-transition space}, and \emph{transition-in-type space} of~$\underline{L}$ are $T\coloneqq L^L$, $\cT \coloneqq \LSh^L$, and $\fT\coloneqq\LSh^{\LSh}$, respectively;
note that  $\T\subseteq \cT$ in a natural fashion (upon identifying $L$ with the subset $\tp(L)\subseteq\fL$).
These spaces generally include (ultra)transitions that are not shard-to-shard.

We regard $\cT$ as the topological product $\prod_{v\in L}\LSh$;
equivalently, via the inclusion $\LSh\subseteq\RR^{\cP}$, the space $\cT$ is topologized as a subspace of the product $\RR^{\cP\times L}$:
This is the topology of pointwise convergence of the real functions $v\mapsto P(f(v))$ for fixed $P\in\cP$.
The space $T\subseteq\cT$ inherits the subspace topology (of pointwise convergence).

\subsubsection{Extendable Layer Transitions}
\label{sec:Sh-extend-comps}
The Continuity Axiom ensures that every computation $\gamma(\cdot): L\to L$ is continuous;
however, it need not extend to a continuous map $\LSh\to\LSh$,
which renders such realized computations rather poor foundational blocks for our subsequent treatment of deep computations.
To remedy such deficiency, we will axiomatically require that realized computations be extendable as suggested in Remark~\ref{rem:non-conts-extn}.
Such an extendability requirement is rather strong;
moreover, its consequences are strongest in regard to the restrictions to compacta of t-ts, rather than the t-ts themselves.
In this light, it is natural to require computations to only be extendable to compact shards~$\fL[\rb]$ when restricted to shards~$L[\rb]$, which motivates the following definition.

A transition-in-type $\ff: \LSh\to\LSh$ is called \emph{shard-continuous} (or \emph{Sh-continuous}) if its restriction to each shard $\fL[\rb]$ is a continuous map $\fL[\rb]\to\fL[\sbb]$ into some type-shard $\fL[\sbb]$
(in particular, a Sh-continuous t-t is shard-to-shard).
A transition $f$ (more generally, a u-t $f$) is called \emph{Sh-extendable} if, for every sizer~$\rb$, its restriction $f{\restriction}L[\rb]$ extends to a continuous function $\fL[\rb]\to\fL[\sbb]$ into some type-shard.
(It suffices to impose this condition for sizers~$\rb$ in a given exhaustive collection~$R$.)

\begin{remark}\label{rem:Sh-ext-to-conts}
  A continuous shard-to-shard t-t is Sh-continuous, but the converse fails in general.

\end{remark}

\subsubsection{Spaces of transitions-in-type}\label{sec:spac-trans-type}

Both $T$ and $\fT$ are semigroups under the binary operation $(f,g)\mapsto f\circ g$ of composition;
this operation is continuous in the left argument~$f$, but \emph{not} in the right argument~$g$.

The subsets $\TShT\subseteq T$ (resp., $\TSh\subseteq\fT$) of transitions (resp., t-ts) that are shard-to-shard are sub-semigroups;
however, $\TShT,\TSh$ are typically \emph{not} closed subspaces.

Recall that a set $R$ of sizers is \emph{exhaustive} if,
for any sizer~$\sbb$, there exists $\rb\in R$ such that $r_P\ge s_P$ for all $P\in\cP$ (cf., Remark~\ref{rem:cofinal-sizer-set}).
In particular, $\LSh = \bigcup_{\rb\in R}\fL[\rb]$ in such case.%
\footnote{\emph{Provided $\cP$ is countable}, by Proposition~\ref{prop:standard-types-countable-P}, we have $\fL = \LSh = \bigcup_{\rb\in R}\fL[\rb]$ for $R$ exhaustive.}

Given a sizer~$\rb$, we say that $\ff\in\fT$ is \emph{$\rb$-preserving} if it restricts to a map $\fL[\rb]\to \fL[\rb]$.
The set of $\rb$-preserving transitions $\ff\in\fT$ is denoted $\fT[\rb]$.
A collection $F\subseteq\fT[\rb]$ is called \emph{$\rb$-preserving}.

Given an exhaustive collection~$R$ of sizers, we say that \emph{$R$ confines~$\ff\in\fT$}, or \emph{$\ff$ is~$R$-confined}, if $\ff$ is $\rb$-preserving for each $\rb\in R$.
The collection of all $R$-confined t-ts is denoted~$\fT[R]$;
it is a closed sub-semigroup of~$\fT$.
One sees that $\fT[R]\subseteq\TSh$, i.e., $R$-confined t-ts are necessarily shard-to-shard.
Moreover, $\fT[R]$ is compact as shown in Proposition~\ref{prop:sh2sh-compactness} below.

The notions above have formally identical analogues for transitions, i.e., with $T, L, L[\rb]$ in place of~$\fT, \LSh, \fL[\rb]$.

 A family $F\subseteq\fT$ is:
\begin{itemize}
\item \emph{confined by an exhaustive sizer collection~$R$} (or: \emph{$R$-confined}) if~$F\subseteq\fT[R]$;
\item \emph{pointwise bounded at $\fv\in\LSh$} if there is a sizer~$\sbb$ (a \emph{pointwise bound} for $F$ at~$\fv$) such that $\ff(\fv)\in\fL[\sbb]$ for all $\ff\in F$;
\item \emph{pointwise bounded on $S\subseteq\LSh$}, if it is pointwise bounded at every~$\fv\in S$;
\item \emph{pointwise bounded}, if it is pointwise bounded on~$\LSh$.
\end{itemize}

\begin{remarks}\label{rem:pw-bdd-t-t-sets}
  \begin{enumerate}
  \item Given $R$ exhaustive, we have $\fT[R] \subseteq \fT[\sbb^{(\cdot)}]$ for the pointwise bounds $\sbb^{(\cdot)} = (\sbb^{(\fv)})_{\fv\in\LSh}$ given by
\begin{equation}\label{eq:pw-bds-from-R}
  \sbb^{(\fv)} \coloneqq \rb^{(\fv)},
  \qquad\text{where $\rb^{(\fv)}\in R$ is any sizer such that $\fv\in\fL[\rb^{(\fv)}]$}
\end{equation}
(at least one such sizer exists, since $\LSh = \bigcup_{\rb\in R}\fL[\rb]$ for exhaustive~$R$).%
  \footnote{Such pointwise bounds are \emph{not} canonically determined by~$R$;
    in particular, one may not instead let $s^{(\fv)}_P \coloneqq \inf\{r_P: \text{$\rb\in R$ and $\fv\in\fL[\rb]$}\}$, for such infimum need not be attained by any sizer of~$R$, and the type-shard $\fL[\sbb]$ of the resulting sizer~$\sbb$ may be strictly smaller than $\fL\cap\RR{[\sbb]}$ (cf., §~\ref{sec:shards})
    —it may even be empty.}
In particular, every confined family $F\subseteq\fT$ is pointwise bounded.
  \item $F\subseteq\fT$ is pointwise bounded on~$S\subseteq\LSh$ iff there is a collection $\sbb^{(\cdot)} \coloneqq (\sbb^{(\fv)})_{\fv\in S}$ of pointwise bounds at each point~$\fv\in S$.
  The corresponding set of t-ts is denoted $\fT[\sbb^{(\cdot)}]$;
  thus
\begin{equation}\label{eq:pointwise-bounded-trns-spc}
  \fT[\sbb^{(\cdot)}]
  \coloneqq
  \{\ff\in\fT: \text{$\ff(\fv)\in \fL[\sbb^{(\fv)}]$ for all $\fv\in S$}\}\quad
  \left( = \prod_{\fv\in S}\fL[\sbb^{(\fv)}]\times \prod_{\fv\in\LSh\setminus S}\LSh \right).
\end{equation}
\end{enumerate}
\end{remarks}

The notions of \emph{$\rb$-preserving}, \emph{$R$-confined}, and \emph{pointwise bounded} (shard-to-shard) layer transitions (and collections of such transitions), as well as the transition spaces $T[\rb]$, $T[R]$, $T[\sbb^{(\cdot)}]\subseteq T$, are obtained from those for transitions-in-type \emph{mutatis mutandis}
(simply replacing $\LSh$ by~$L$, and $\fT$ by~$T$).

\begin{proposition}\label{prop:sh2sh-compactness}
  For any collection $\sbb^{(\cdot)} = (\sbb^{(\fv)})_{\fv\in\LSh}$ of sizers at \emph{all} points $\fv\in\LSh$, the space $\fT[\sbb^{(\cdot)}]$ of $\sbb^{(\cdot)}$-pointwise bounded transitions-in-type is compact.
  In particular, $\fT[R]$ is compact for any exhaustive sizer collection~$R$.
\end{proposition}
\begin{proof}
  The collection $(\sbb^{(\cdot)})$ specifies pointwise bounds at all points $\fv\in\LSh$, hence the product space in~\eqref{eq:pointwise-bounded-trns-spc} above is compact, by Tychonoff’s Theorem (being a product of compact factors $\fL[\sbb^{(\fv)}]$ only).

  If $R$ is exhaustive, we have $\fT[R] \subseteq \fT[\sbb^{(\cdot)}]$ for pointwise bounds $\sbb^{(\cdot)}$ as in~\eqref{eq:pw-bds-from-R};
  since $\fT[R]$ is closed in~$\fT$, it is therefore a closed subset of the compact space~$\fT[\sbb^{(\cdot)}]$, hence compact.
\end{proof}

\subsection{Computations and ultracomputations (deep computations)}
\label{sec:computation-types}

\subsubsection{The Extendability Axiom}
\label{sec:Ax-Extend}

The \emph{transition $\ev(\gamma,\cdot)$ associated to~$\gamma\in\Gamma$} —also denoted $\gamma(\cdot)$— of a layer transformation $\gamma\in\Gamma$ is the map $v\mapsto\ev(\gamma,v)$.

If such a transition~$\gamma(\cdot)$ is Sh-extendable, we call it the \emph{computation by~$\gamma$}, or \emph{realized by~$\gamma$} for emphasis.
By an abuse of notation, we will denote the extension $\LSh\to\LSh$ still by $\gamma(\cdot)\in \fT$.

For the remainder of this paper, we assume that CCSs satisfy the following
\begin{itemize}
\item \textbf{Extendability Axiom.}\quad
  \emph{Each layer transformation $\gamma\in\Gamma$ induces a Sh-extendable computation~$\gamma(\cdot)$.}
\end{itemize}

The Extendability Axiom gives a natural (injective) map $\Gamma\to\fT: \gamma\mapsto \gamma(\cdot)$.
The semigroup~$\Gamma$ is topologized via (the pullback of) this map, i.e., by the topology of pointwise convergence;
it is the subspace topology obtained upon identifying $\Gamma$ with the set $\Gamma(\cdot) \coloneqq \{\gamma(\cdot): \gamma\in\Gamma\} \subseteq \fT$, called the \emph{space of realized computations}.

It follows from the Reduction Axiom that the above topology on~$\Gamma$ is Hausdorff.

\subsubsection{Realized vs.\ deep computations}\label{sec:deep-computations}

The \emph{space $\fD$ of ultracomputations} is the topological closure $\overline{\Gamma(\cdot)} \subseteq \fT$.
A transition $\ff\in\fD$ will be called a \emph{deep computation}, \emph{ultracomputation}, or \emph{ucomp} for short.
Although any computation is a deep computation in its own right, the adjective “deep” implies that $\ff$ may be an \emph{un}realized computation, i.e., \emph{not} of the form~$\gamma(\cdot)$.
Deep computations are typically (Sh-)\emph{dis}continuous layer transitions.
Even if Sh-continuous, an ultracomputation may be unrealized.

Every deep computation is of the form $\ff_{\cU}\coloneqq \Ulim_i\gamma_i: v\mapsto \ff_{\cU}(v) \coloneqq \Ulim_i\gamma_i(v)$ for some indexed family $\gammab\coloneqq (\gamma_i)_{i\in I}\subseteq\Gamma$ and some ultrafilter $\cU$ on~$I$.
(Without loss of generality, one may always take $\cU$ as an ultrafilter on $\Gamma$ itself.)%
\footnote{For arbitrary $\cU$ on (say)~$\Gamma$, the ultracomputation $\ff_{\cU}$ need not be defined:
$(\gamma(v))$ might not $\cU$-converge for certain~$v\in L$.}

For any sizer collection $\sbb^{(\cdot)}$, let $\fD[\sbb^{(\cdot)}] \coloneqq \fT[\sbb^{(\cdot)}]\cap\fD$ be the set of ultracomputations with pointwise bounds~$\sbb^{(\cdot)}$.
Since $\fD\subseteq\fT$ is closed by definition, the space $\fD[\sbb^{(\cdot)}]$ is also closed in~$\fT$.
For any fixed sizer $\rb$ and exhaustive~$R$, we see that $\fD[\rb] \coloneqq \fT[\rb]\cap\fD$ and $\fD[R] \coloneqq \fT[R]\cap\fD$ (the sets of ultracomputations $\rb$-preserving and $R$-confined, respectively) are also closed.

By an abuse of nomenclature, we say that an element $\gamma\in\Gamma$ \emph{admits pointwise bounds~$\sbb^{(\cdot)}$} (resp., \emph{is $\rb$-preserving}, \emph{is $R$-confined}) if its transition type $\gamma(\cdot)\in T$ does (resp., is).
We denote by $\Gamma[\sbb^{(\cdot)}]$, $\Gamma[\rb]$, and~$\Gamma[R]$, respectively, the sets of transformations~$\gamma\in\Gamma$ with associated transitions in $\fT[\sbb^{(\cdot)}]$, $\fT[\rb]$, and $\fT[R]$.
The respective uniform notions as $\gamma$ varies in some subset~$\Delta\subseteq\Gamma$ become:
\emph{$\Delta$ admits uniform pointwise bounds~$\sbb^{(\cdot)}$}, \emph{is $\rb$-preserving}, and \emph{is $R$-confined}, respectively.

An ultracomputation $\ff:L\to L$ with values in $L\subseteq\fL$ is called \emph{quasi-realized};
these constitute the set $\cD = \fD\cap T$:
the space of \emph{quasi-realized ultracomputations}.
Let $\cD[\sbb^{(\cdot)}]\coloneqq \fD[\sbb^{(\cdot)}]\cap T$, $\cD[\rb]\coloneqq \fD[\rb]\cap T$, and $\cD[R]\coloneqq \fD[R]\cap T$.

\begin{proposition}\label{prop:transf-type-spaces}
  For any sizer $\rb$ and exhaustive collection~$R$:
  \begin{enumerate}
  \item each of the sets $\Gamma[\rb]$, $\Gamma[R]$ is a sub-semigroup of~$\Gamma$, and is a closed subset of~$\Gamma$;
  \item $\fD[\rb]$, $\fD$ are closed sub-semigroups of~$\fT$;
  \item $\fD[R]$ is a compact sub-semigroup of~$\fT$;
  \item $\cD[\rb]$, $\cD[R]$, $\cD$ are closed sub-semigroups of~$T$.
  \end{enumerate}
\end{proposition}
The ultracomputation space $\fD$ is akin to the concept of “enveloping group” (of~$\Gamma(\cdot)\subseteq\fT$).
However, only the confined sub-semigroups $\fD[R]$ are compact (the full space~$\fD$ is typically noncompact).
\begin{proof}
  The set $\beta\Gamma$ of ultrafilters on~$\Gamma$ is itself a semigroup under a natural (“convolution”) operation $(\cU,\cV)\mapsto\cU{*}\cV$~\cite{HindmanStrauss2010}.
  This operation of convolution possesses (and is essentially characterized by) the following property
 —when $\Gamma$ is identified with the transitions semigroup~$\Gamma(\cdot)$:
  If two transitions are of the form $\ff_{\cU}: \fv\mapsto\Ulim_{\gamma}\gamma(\fv)$, $\ff_{\cV}: \fv\mapsto\Vlim_{\gamma}\gamma(\fv)$, then $\ff_{\cU}\circ \ff_{\cV} = \ff_{\cU{*}\cV}$.
  It follows that $\fD\subseteq \fT$ is a sub-semigroup.
  As the intersection of a compact subset (by Proposition~\ref{prop:sh2sh-compactness}) with a closed subset of~$\fT$, we see that $\fD[R] = \fT[R]\cap\fD$ is compact.
  The remaining topological statements are all trivial and left to the reader.
\end{proof}

\section{Deep iterations and deep equilibria}
\label{S:DEQ}

Throughout this section, fix a CCS $\cC = \langle\underline{L}, \underline{\Gamma}\rangle$.
For convenience, we assume some element~$\id\in\Gamma$ (“identity”) satisfies the equality $\id(v) = v$ for all $v\in L$.

We reiterate the Extendability Axiom that each layer transition $\gamma\in\Gamma$ extends to a Sh-continuous transition $\gamma(\cdot)\in\TSh$.

\subsection{Layered and iterative computations}
\label{sec:LCs-ICs}
Let $\gammab = (\gamma_n)_{n\in\NN} \subseteq \Gamma$ be any sequence of computations
(i.e., any element of the product space $\Gamma^{\omega}\coloneqq \prod_{n\in\NN}\Gamma$).
We regard $\gammab$ as a sequence of “computation steps” to be successively applied (see the definition of Layered Computation below).
The computation $\gamma_n$ will be called the \emph{$n$-th atomic step}, or the \emph{transition at layer~$n$ (to layer~$n+1$)}.

\paragraph{Layered Computations (LCs)}
Given a sequence $\gammab\in\Gamma^{\omega}$ of computation steps, the sequence $\gammab^{(\circ)} = (\gammab^{(n)})_{n\in\NN} \in \Gamma^{\omega}$ defined recursively by
\begin{align*}
  \gammab^{(0)} &\coloneqq \id, \\
  \gammab^{(n+1)} &\coloneqq \gamma_n\gammab^{(n)}\qquad\text{for all $n\in\NN$,}
\end{align*}
(i.e., $\gammab^{(n)} \coloneqq \gamma_{n-1}\gamma_{n-2}\dots \gamma_1\gamma_0$)
is called the \emph{layered computation with atomic steps~$\gammab$} (or \emph{LC$\gammab$}, for short).%
\footnote{Thus, LC$\gammab$ denotes $\gamma^{(\circ)}$, simply adding context to indicate the layer transitions~$\gammab$ yielding $\gamma^{(\circ)}$.}
The term $\gammab^{(n)}$ is called the \emph{$n$-composite computation step} of~LC$\gammab$.
The set of layered computations LC$\gammab$ obtained as $\gammab\in\Gamma^{\omega}$ varies is denoted~$\Gamma^{(\circ)}$.

For a sizer~$\rb$, let $\Gamma^{(\circ)}_{[\rb]} \coloneqq \Gamma^{(\circ)}\cap(\Gamma[\rb])^{\omega}$ be the set of $\rb$-preserving LCs.
(Note that it is the products $\gammab^{(n)}$ —but not necessarily the atomic steps~$\gamma_n$— that are required to preserve the layer~$L[\rb]$.)
For an exhaustive sizer family~$R$, let  $\Gamma^{(\circ)}_{[R]} \coloneqq \Gamma^{(\circ)}\cap(\Gamma[R])^{\omega}$ be the set of \emph{$R$-confined LCs} (or \emph{LCs confined by~$R$}).

The \emph{LC$\gammab$-evolution of a state~$v\in L$} is the sequence
\begin{equation*}
  \gamma^{(\circ)}(v) \coloneqq (\gammab^{(n)}(v))_{n\in\NN} = (v, \gamma_0(v), \gamma_1\gamma_0(v), \gamma_2\gamma_1\gamma_0(v), \dots).
\end{equation*}
The term “evolution” means “$\gammab$-evolution” henceforth, whenever $\gammab$ is given by context.
The \emph{state at stage~$n$ of~$v$} under evolution is $\ev(\gamma^{(n)},v)$.

\paragraph{Iterative computations (ICs)}
Any fixed $\gamma\in\Gamma$ yields a constant sequence $\gammab = (\gamma)_{n\in\NN}$.
The corresponding LC has composite steps given by the sequence $(\gamma^n)_{n\in\NN}$ of compositional powers (iterates) of~$\gamma$;
we will call such LC an \emph{iterative computation \emph{(or just \emph{iteration})} by~$\gamma$}, and denote it by IC$\gamma$.
It is appropriate to think of iterative computations as evolving by “tying parameters” in the sense that all atomic steps are always the same~$\gamma$ (i.e., the “tied parameter” is~$\gamma$ itself).
Note that IC$\gamma$ is $R$-confined (or $\rb$-preserving) if and only if the fixed atomic step~$\gamma$ is so.

\subsection{Deep layers, deep iterates, and equilibria}
\label{sec:deep-layers}

\subsubsection{Deep layers}
A \emph{deep layer} of LC$\gammab$ is any deep computation that is an accumulation point of the sequence of composites $(\gammab^{(n)}(\cdot))_{n\in\NN} \subseteq \fT$.
Any such deep layer is of the form $\gammab^{(\cU)}: v\mapsto\Ulim_n\gamma^{(n)}(v)$
obtained as (pointwise) $\cU$-ultralimit via a nonprincipal ultrafilter $\cU$ on~$\NN$.
(We use the notation $\gammab^{(\cU)}$ for such ultracomputation (in-type) when the dependence on $\gammab$ and~$\cU$ is to be made explicit.)
For such a confined LC$\gammab$, deep limits exist for arbitrary~$\cU$, by Proposition~\ref{prop:existence-ucomps}.
If LC$\gammab$ is not confined, the computations $\gamma^{(n)}$ may diverge.

\subsubsection{Deep iterates}
\label{sec:deep-iterates}
A \emph{deep iterate} of $\gamma\in\Gamma$ is a deep layer for IC$\gamma$.

The deep layer that is obtained via an ultrafilter $\cU$ on~$\NN$ is denoted~$\gamma^{(\cU)}$;
it need not exist in general, but does if~$\gamma$ is confined (by Proposition~\ref{prop:existence-ucomps}).
Every deep iterate is a deep computation.

\begin{remark}\label{rem:disconts-ufns}
In the nomenclature of \cite{BaiKolterKoltun2020}, a deep layer of LC$\gammab$ is an “implicit layer”.
They consider primarily compositions of layer transitions (i.e., LCs in our sense) with “tied parameter”~$\gamma$ (the same layer transition at each stage), i.e., ICs in our sense.
From our perspective, implicit layers are given each by some nonprincipal ultrafilter~$\cU$ on~$\NN$, i.e., are of the form $\gammab^{(\cU)}$.
\end{remark}

\subsubsection{Deep equilibria}
\label{sec:deep-equilibria}

A \emph{deep equilibrium (layer)} of IC$\gamma$ is an idempotent deep iterate $\fid = \gamma^{(\cU)}\in\fD$ ($\subseteq \fT$), i.e., a deep iterate~$\fid$ such that $\fid(\fid(\fv)) = \fid(\fv)$ for all $\fv\in \LSh$ (hence the nomenclature “equilibrium”).
It will also be called a \emph{(deep) iterative equilibrium} of~$\gamma$.

\begin{remark}\label{rem:DEQ-self-reference}
 Although any iterative equilibrium~$\fid$ of~IC$\gamma$ satisfies $\fid\circ\fid = \fid$, one generally has $\gamma\circ\fid\ne\fid\ne\fid\circ\gamma$.
  The “equilibrium” property is \emph{self-referential}, rather than in direct relation to the original computation~$\gamma$.
  Let us call a deep iterate $\gamma^{*}$ of IC$\gamma$ “$\gamma$-fixed” if $\gamma\circ\gamma^{*} = \gamma^{*} = \gamma^{*}\circ\gamma$.
  Such $\gamma$-fixed deep iterates need not exist even under the strong hypothesis (ensuring that deep iterates exist at all) that IC$\gamma$ is confined.
  On the other hand, if perchance a deep iterate $\gamma^{*}$ of IC$\gamma$ satisfies $\gamma\circ \gamma^{*} = \gamma^{*}$, then certainly $\gamma^{*}$ is a deep equilibrium in our sense.
\end{remark}

\begin{therm}[Existence of Deep Iterative Equilibria]
  \label{thm:DEQ-exist}
  Let $\gamma\in\Gamma$ be confined.
  Then there exists at least one iterative equilibrium~$\fid$ for~IC$\gamma$.
\end{therm}
Theorem~\ref{thm:DEQ-exist} is essentially a particular case of the classical Ellis--Numakura Lemma;
the proof below is standard (as in~\cite{Furstenberg:1981}).

One cannot generally hope that deep iterative equilibria exist without some boundedness assumption (such as confinement).
Moreover, $\fid{\restriction}L$ need not take values in~$L$, so it is not even composable with itself \emph{a priori}!
This highlights the need to consider transitions \emph{in type} rather than as maps $L\to L$ on the layer state space~$L$.
\begin{proof}
  Let $R$ confine~IC$\gamma$.
  Let $G \subseteq \fT[R]$ be the topological closure of the semigroup $\{\gamma^n(\cdot): n\ge 1\}\subseteq T$ of transitions by iterates of~$\gamma$ (excluding the trivial iterate $\gamma^0 = \id$).
  By Proposition~\ref{prop:transf-type-spaces}
  (in the CCS obtained from~$\cC$ with computations semigroup $\langle\gamma\rangle$ generated by~$\gamma$),
  $G$ is a compact Hausdorff topological semigroup under composition $(\ff,\fg)\mapsto \ff\circ\fg$, which is continuous in the left argument~$\ff$ (for fixed~$g$).
  Elementary algebraic and topological considerations (in particular, the compactness of~$G$), and Zorn's Lemma, imply that $G$ has some minimal closed (nonempty) sub-semigroup~$H$ (i.e., $H$ has no proper closed sub-semigroups).
  Fix $\fid\in H$.
  The set $H\circ\fid$ is closed (since $g\mapsto g\circ\fid$ is continuous and $H$ is closed);
  moreover, $(H\circ\fid)\circ(H\circ\fid)\subseteq H\circ\fid$, so $H\circ\fid\subseteq H$ is a closed sub-semigroup, hence $H\circ\fid = H$ by minimality of~$H$.
  Therefore, $\ff\circ\fid = \fid$ for some $\ff\in H$.
  Let $H' \coloneqq \{\fg\in H: \fg\circ\fid = \fid\}\ni\ff$.
  Thus, $H'$ is clearly a nonempty sub-semigroup of~$H$, and also closed
  (as the inverse image of the closed singleton $\{\fid\}$ under the continuous map $\fg\mapsto \fg\circ\fid$, again).
  By minimality, $H'=H\ni\fid$, so $\fid\circ\fid = \fid$.
\end{proof}

\subsection{Examples and discussion of deep iterates and deep equilibria}
\label{sec:deep-equi-exa}

\begin{example}\label{exa:DEQ-finite-set}
  Let $L$ be a finite set of, say, $m\ge 1$ distinct elements.
  The choice of predicates is inessential in this context:
  we may simply take $L = [m] \coloneqq \{1,\dots,m\}$:
  it is finite and therefore realcompact.
  Let $f:[m]\to[m]$ be any function, and $\underline{\Gamma} = \langle f\rangle \coloneqq \{f^n:n\in\NN\}$ (as a semigroup under composition) act on $[m]$ by functional application~$\ev:(g,i)\mapsto g(i)$.
  Let $P_{\id}: L\to\RR: i\mapsto i$ (the identity function) be the sole predicate on~$[m]$.
  In this way, we obtain a (realcompact) CCS $\cC = (([m],P_{\id}), \underline{\Gamma})$.
  Since~$[m]$ is finite, there is $n\ge 1$ such that $S\coloneqq f^n([m]) = f^{n+1}([m]) = f(S)$ (thus, $S\ne\emptyset$);
  in particular, $f$ restricts to a bijection of $S\to S$;
  by relabeling points of~$L=[m]$ if necessary, we may as well assume $S = [k]$ ($1\le k\le m$).
  Thus, $g \coloneqq f\restriction[k]$ is a permutation of~$[k]$.
  Let $K$ be the order of~$g$ (thus, $1\le K\le k!$).
  Let $N$ be any integer such that $N\ge n$ and~$K$ divides~$N$.
  Then $f^{*}=f^N$ is a deep iterative equilibrium of~$f$:
  indeed, for $i\in[m]$,
  \begin{equation*}
    \begin{split}
      f^{*}(f^{*}(i)) &
                        = f^N(f^{N-n}(f^n(i)))
                        = g^N(g^{N-n}(f^n(i)))\\
                      &\qquad\text{(since $f^n([m])\subseteq[k]$ and $g = f\restriction[k]$)}\\
                      &= g^{N-n}(g^N(f^n(i)))
                        = g^{N-n}(f^n(i))\\
                      &\qquad\text{(since $K$ divides $N$ and $g^K$ is the identity)}\\
                      &= f^{N-n}(f^n(i))
                        = f^N(i) = f^{*}(i).
    \end{split}
  \end{equation*}
  It is easy to show that $f^{*}$ is the unique iterative equilibrium of~$f$ in such case.
\end{example}

\begin{example}\label{exa:DEQ-unit-interval}
  Consider CCSs of the form $\cC = \langle([0,1], \{P_{\id}\}), \underline{\Gamma}\rangle$ as in~\ref{sec:unit-interval}, where $\gamma: [0,1]\to[0,1]$ is a continuous map, $\Gamma=\langle \gamma\rangle$ the semigroup of iterates of~$\gamma$ under composition, acting by functional application on~$[0,1]$.
  Already in this one-dimensional compact setting, there is a variety of possible behaviors of deep iterates and equilibria of~IC$\gamma$.

  If $\Gamma(\cdot)$ is an equicontinuous family of functions on~$[0,1]$, the Arzelà--Ascoli Theorem implies that there exists a (sub)sequence $(\gamma^{n_k})_{k\in\NN}$ of iterates converging \emph{uniformly} to a \emph{continuous} limit~$\bar{\gamma}: [0,1]\to[0,1]$, which is therefore a continuous deep iterate of~$\gamma$.
  In general, however, even if some deep iterates~$\bar{\gamma}$ are continuous, some deep equilibria may be discontinuous.
  Typically (and necessarily so when $\gamma$ is a chaotic function —e.g., the logistic map $\gamma(v) = 4v(1-v)$), the semigroup $\Gamma(\cdot)$ is not an equicontinuous collection of functions, and deep equilibria (as well as deep iterates) are necessarily discontinuous.
  Moreover (in contrast with the equicontinuous case possessing continuous deep iterates achieved sequentially),
  deep iterates~$\ff$ of a chaotic~IC$\gamma$ \emph{cannot} be obtained as sequential limits $\lim_k\gamma^{n_k}$, but generally only as \emph{ultra}limits.
\end{example}

\begin{figure}[htbp]
  
  \centering
  
    \includegraphics[width=0.55\textwidth]{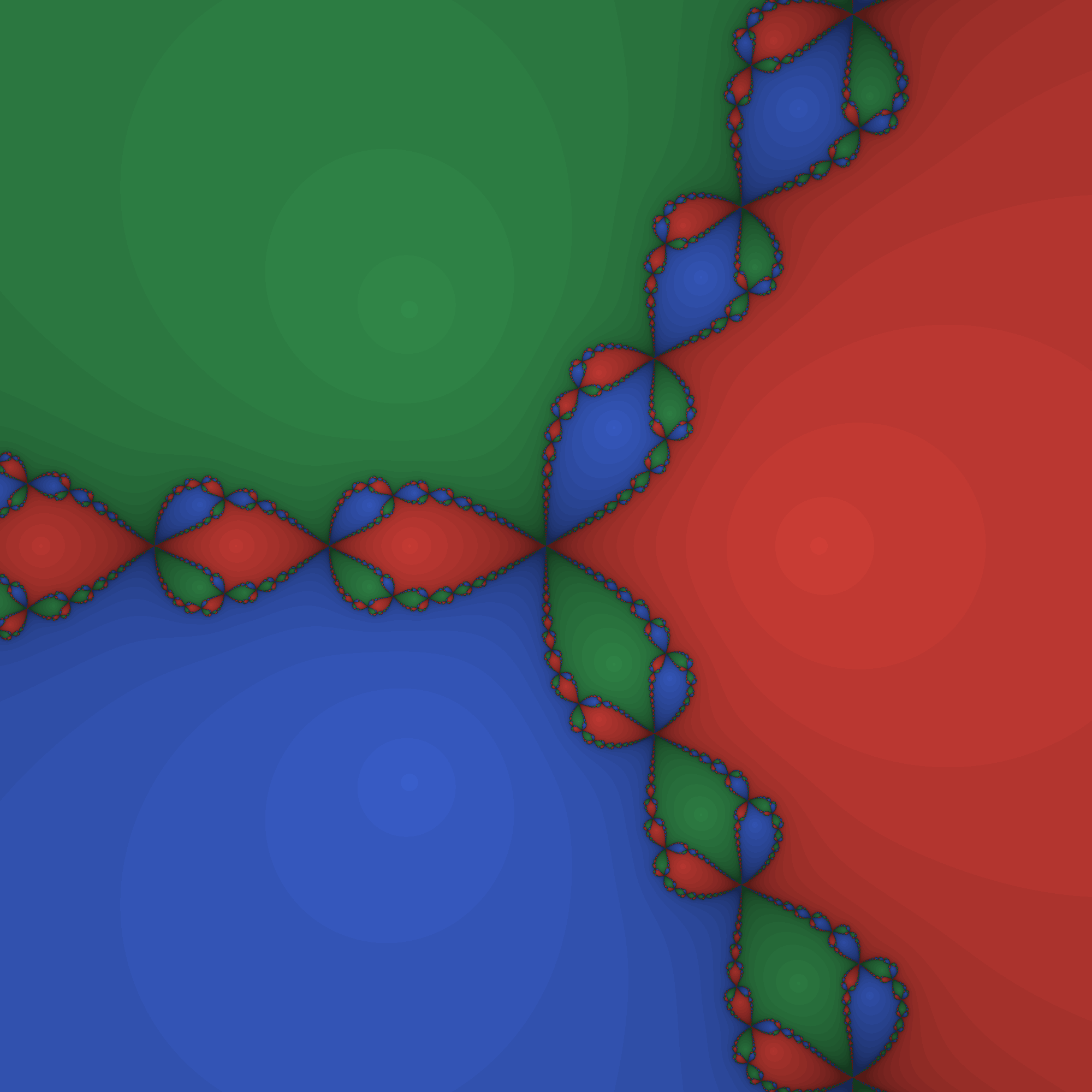}
  
  \caption{The Newton fractal of $p(z) = z^3-1$:
    the three basins of attraction of the cube roots of unity $1$ and $\zeta_{\pm} = \exp(\pm2\pi i/3)$ are shown  in colors (red for~$1$, green  for~$\zeta_+$, blue for~$\zeta_-$);
    the intricate boundary $K$ separating the basins is a Newton fractal (Example~\ref{exa:DEQ-Newton}).
    It is impossible to give an effectively computable description of $\gamma^{*}$ on~$K$;
    its proof of existence is inexplicit, requiring the abstract formalism of ultralimits.}
  \label{fig:newton-fractal}
\end{figure}

\begin{example}[Deep iterates and equilibria of Newton’s Method] \label{exa:DEQ-Newton}
  Fix a polynomial $p$ with (real or) complex coefficients and degree $d\ge 2$.
  Consider the CCS
  \begin{equation*}
  \bigl\langle(\hat{\CC}, \{U,V,W\}), \langle \gamma\rangle\bigr\rangle,
\end{equation*}
where
\begin{itemize}
  \item $\hat{\CC} = \CC\cup\{\infty\}$ is the Riemann sphere, which we identify with the unit sphere $S^2=\{(u,v,w)\in\RR^3: u^2+v^2+w^2=1\}$ via, e.g., the stereographic projection $(u,v,w)\mapsto z = (u+iv)/(1-w)$ (and $(0,0,1)\mapsto\infty$);
  \item $\hat{\CC}\to S^2: z\mapsto(U,V,W)$ is the inverse of the stereographic projection, regarded as a triple of predicates $U,V,W: \hat{\CC}\to[-1,1]$; and
  \item
    \begin{equation*}
      \gamma(z)
      \coloneqq
      \begin{cases}
        z - \frac{p(z)}{p'(z)} & \text{($p'(z)\ne 0$)}\\
        z & \text{($p'(z)=0=p(z)$)}\\
        \infty & \text{($p'(z)= 0 \ne p(z)$, or $z=\infty$)} %
      \end{cases}
    \end{equation*}
    is the transition (“Newton map”) carrying out one step of Newton’s method to find the roots of~$p(z)$, extended to a continuous map of~$\hat{\CC}$ by $p(\infty) \coloneqq \infty$.
    Under the assumption $\deg p\ge 2$, the map $z\mapsto p(z)/p'(z)$ also extends continuously to~$\hat{\CC}$, by $z\mapsto\infty$ when either $p'(z) = 0\ne p(z)$ or $z=\infty$, and by $z\mapsto 0$ at a multiple root of~$p$ (where $p'(z) = 0 = p(z)$).
  \end{itemize}
  Since $\hat{\CC}$ is compact and $z\mapsto(U(z),V(z),W(z))$ is a homeomorphic embedding, $\hat{\CC}$ is realcompact.
  In fact, $\hat{\CC}$ is equal to the shard $\hat{\CC}[1,1,1] = \bigl\{z: \max\bigl(\nrm{U(z)},\nrm{V(z)},\nrm{W(z)}\bigr) \le 1\bigr\}$;
  thus, $\gamma$ is automatically confined (by $R$ consisting of the single sizer $\rb = (r_U,r_V,r_W) = (1,1,1)$).%
  \footnote{Most of the discussion of deep iterates $\gamma^{*}$ below remains valid upon substituting a rational complex map for~$p$
    —except that the point~$\infty$ is effectively stripped of any special role.}

  Let $\gamma^{*}$ be any deep iterate of~$\gamma$.
  At any point $z\in\hat{\CC}$ for which the Newton method converges to a root $w$ of~$p(z)$ (in particular, at any $z$ sufficiently close to a simple such root~$w$), we have $\gamma^{*}(z) = w$ ($=\gamma(w)$, since $p(w)=0$)
  —in fact, $\gamma^{*}$ is \emph{constant} (with value the same~$w$) in some neighborhood of~$z$.
  We also have $\gamma^{*}(\infty) = \infty = \gamma(\infty)$;
  however, $\infty$ is a repeller
  (this follows from an easy calculation showing that $\gamma(z)/z\to 1-d^{-1}$ as $\nrm{z}\to+\infty$).%
  \footnote{The fixed repeller~$\infty$ is typically the (terminal) point of countably many \emph{finite} orbits $(\gamma^n(z))_{n\in\NN}\subseteq\CC$ ending there.}
  In general, however, $w\coloneqq \gamma^{*}(z)$ is not necessarily a root of~$p$. %

  The Fatou set of~$\gamma$ consists of those points having a neighborhood on which some deep equilibrium is attained as a uniform sub-sequential limit of~$(\gamma^n)$
  —then so is the case for every deep equilibrium;
  its complement is the Julia set, consisting of points in whose every neighborhood no subsequence of $(\gamma^n)$ converges uniformly
  (or converges pointwise at all, it turns out).
  Cf.~\cite[§~4]{Milnor:2006} for this standard characterization of the Fatou and Julia sets.
  The Julia set of the Newton map~$\gamma$ is called the \emph{Newton fractal} of~$p$.

  For the polynomial $p(z) = z^3-1$, the Fatou set is easy to characterize:
  it is the union of the basins of attraction of its three roots $1$ and $\zeta_{\pm} = \exp(\pm2\pi i/3)$ under iterates of the Newton map~$\gamma(z) = (2z + z^{-2})/3$;
  its complement $K$ is the boundary of each and all of these three regions: the (closed nowhere dense) Newton fractal of~$z^{3}-1$ (Figure~\ref{fig:newton-fractal}).
  All deep iterates and equilibria have the same value $\lim_{n\to\infty}\gamma^n(z)$ at (convergent) inputs~$z\in\CC\setminus K$
  —the common restriction of all such $\gamma^{*}$ to $\CC\setminus K$ is continuous;
  however, each deep iterate restricts to an extremely complicated map of $K$ to itself.
\end{example}

\begin{example}\label{exa:DEQ-a-la-BKK}
  The definitions of deep layer state and deep iterative equilibrium above are motivated by the notions of “Deep Equilibrium (DEQ)” in \cite{BaiKolterKoltun2020}.
  However, iterative computations in~\cite{BaiKolterKoltun2020} allow “feeding” the initial state~$v$ as an argument at each application of a (“parameter-tied”, i.e., fixed) layer transformation.
  Capturing deep iterative equilibria in this sense requires generalizing the notion of CCS.
  One way to capture the deep equilibria of Bai \emph{et al.}\ is allowing $\cC$ to be a CCS with $n$-ary (in fact, just binary) layer transformations as in Remark~\ref{rem:n-ary-transfns}(2).
  Indeed, fix a \emph{binary} $\gamma\in\Gamma$, which induces a two-argument layer transition $\gamma(\cdot,\cdot): L\times L\to L$.
  Consider the map $\delta: L\times L\to L\times L$ given by $\delta(v,w)\coloneqq (v,\gamma(v,w))$
  (the first entry of $\delta(v,w)$ is simply a pass-through of the first argument, while the second entry applies the computation~$\gamma$).
  Then the second (nontrivial) entry $w_n \eqqcolon f_n(v)$ of the iterates $\delta^n(v,v) = (v,w_n)$ for $n\in\NN$ represents the evolution of the computation $\gamma$ passing through, at each step, the original argument~$v$ as the first of two inputs.

  If $L$ is realcompact and $\gamma$ is $R$-confined (i.e., restricts to a map $L[\rb]\times L[\rb]\to L[\rb]$ for all $\rb\in R$), the proof of Theorem~\ref{thm:DEQ-exist} is adapted \emph{mutatis mutandis} to computations in CCSs with $n$-ary transitions.
  One shows thus the existence of deep equilibria, i.e., of idempotent maps $\fid: L\to L$ arising as ultralimits of the sequence of iterates~$(f_n)_{n\in\NN}$ of evolution by~$\gamma$.
  (Without a realcompactness assumption, one needs suitable hypotheses on~$\gamma$ akin to Sh-extendability.)

  As an alternative to the use of CCSs with $n$-ary computations, in Appendix~\ref{sec:PFC}, we introduce the notion of \emph{Parametrized Family of Computations (PFC)} to capture computations with feed-through in our framework.
  The ability to compute deep equilibria in an effective sense, as in Bai~\emph{et al.}, presupposes that such equilibria are definable not merely in a continuous, but in a differentiable sense
  (allowing the use of generic solver —or fixed-point— algorithms, which typically rely on gradient-descent methods, e.g., Newton’s algorithm and refinements);
  we explain how such considerations of differentiability may be handled in a CSS with \emph{finitely} many predicates
  (considerations of differentiability when infinitely many observables are involved entail delicate analysis beyond the scope of this paper).
\end{example}

\begin{remarks}\label{rem:disconts-ucomps}
  \begin{enumerate}
  \item The results in §~\ref{sec:DI-computability} below say nothing about effectively computable features %
    of (shard-)\emph{dis}continuous deep iterates or equilibria such as those arising from Newton's method iterations in Example~\ref{exa:DEQ-Newton} above.
    In an upcoming article, we extend the present results to discontinuous ultracomputations that are nevertheless \emph{de facto} effectively computable in a localized sense.

  \item Even in situations where, say, a deep iterate does not quite exist, an ultracomputation may have “meaningful deep features” in a sense that we now explain.
  Consider any CCS $\cC = (L, \langle f\rangle, \circ, \cP)$ (not necessarily realcompact), where $f:L\to L$ is any given (continuous) computation.
  For a fixed $Q\in\cP$, say that~$f$ \emph{has uniformly $Q$-bounded iterates on~$v\in L$} if there exists $s = s^{(v)} > 0$ such that $\nrm{Q(f^n(v))}\le s$ for all $n\in\NN$.
  (Note that this hypothesis does not —at all— impose bounds on other entries $P\circ f^n(v)$ for $Q\ne P\in\cP$.)
  If $\cU$ is any nonprincipal ultrafilter on~$\NN$, the iterate boundedness hypothesis and the compactness of intervals $[-s,s]$ imply that $\Ulim Q(f^n(v))$ exists for all $v\in L$.
  In principle, however, the iterates $f^n$ need not $\cU$-converge in $T$ (i.e., pointwise on~$L$) even if $L$ is realcompact, since (for fixed~$v$) the sequence $(f^n(v))_{n\in\NN}$ may not be entry-wise bounded (only bounded “in $Q$-th entry”, so to speak).
\end{enumerate}
\end{remarks}

  The study of aspects of deep equilibria introduced in Remarks~\ref{rem:disconts-ucomps} is quite delicate, and exceeds the scope of the present paper.

\begin{proposition}\label{prop:existence-ucomps}
  The ultracomputation $\ff_{\cU}\coloneqq\Ulim_{\gamma}\gamma(\cdot)\in \fT[R]$ exists for any exhaustive~$R$ and any ultrafilter~$\cU$ on $\Gamma[R]$.
\end{proposition}
\begin{proof}
  This is an immediate corollary of Proposition~\ref{prop:transf-type-spaces}.
\end{proof}

\section{Explicit computability}
\label{sec:DI-computability}

Throughout this section, we fix a CCS $\cC = \langle \underline{L},\underline{\Gamma}\rangle$.
We assume $\Gamma$ has an identity element $\id$ acting as the identity map on~$L$.
We reiterate the Extendability Axiom that each layer transition $\gamma\in\Gamma$ extends to a Sh-continuous transition $\gamma(\cdot)\in\TSh$.

We shall implicitly identify a predicate symbol~$P$ with the real-valued function $P(\cdot): L\to\RR$ interpreting it in~$\cC$, and also implicitly extend $P(\cdot)$ to a (unique continuous) function $\LSh\to\RR$.

A real function $\varphi: \LSh\to\RR$ will be called \emph{shard-bounded (sh-bdd)} (resp., \emph{Sh-continuous}) if its restriction to each shard $\fL[\rb]$ is bounded (resp., continuous).
(A Sh-continuous such function is necessarily sh-bdd.)

\subsection{Polynomials in predicates and definability. Features of layer transitions.}
\label{sec:definable-predicates}

\begin{itemize}
\item Any predicate $P$ will also be called a \emph{monomial}.%
  \footnote{In real-valued logic, the monomials above are called “atomic”.}
\item A \emph{polynomial} is any function $L\to\RR$ obtained by combining real constants $r\in\RR$ and monomials using any (recursive) combination of the following operations, called \emph{connectives}:
  \begin{itemize}
  \item \emph{Addition:}\quad
    $(\varphi,\psi)\mapsto \varphi+\psi$\quad (where $\varphi+\psi: v\mapsto \varphi(v)+\psi(v)$);
  \item \emph{Multiplication:}\quad
  $(\varphi,\psi)\mapsto \varphi\psi$\quad (where $\varphi\psi: v\mapsto \varphi(v)\psi(v)$).
  \end{itemize}
  The monomials appearing in an expression of some polynomial may be called its \emph{atoms}.%
  \footnote{A polynomial~$\varphi$ need not have a unique expression in terms of monomials, so it is more accurate to say that $\varphi$ has an expression involving certain specific monomials.\label{fn:param-set}}

\item A \emph{definable predicate} is a function $\varphi: L\to\RR$ whose restriction to an arbitrary shard $L[\rb]$ is uniformly approximable by polynomials;%
  \footnote{The notion of “definable predicate” above is less restrictive than the (most) standard one in real-valued logic, wherein approximability is required to hold uniformly over the entire set (“universe”)~$L$.}
  thus, $\varphi$ is definable iff for every $\varepsilon>0$ and sizer $\rb$ there exists a polynomial~$\psi = \psi_{[\rb]}^{(\varepsilon)}$ such that $\nrm{\varphi(v)-\psi(v)}<\varepsilon$ for all~$v\in L[\rb]$.
  The family $(\psi_{[\rb]}^{(\varepsilon)}: \varepsilon>0, \text{$\rb$ sizer})$ is called a \emph{definition scheme for~$\varphi$}.

\end{itemize}
We \emph{only} require definable predicates to be uniformly approximable \emph{on shards} —not uniformly on the full state space~$L$.

\begin{remarks}\label{rem:definability}
  \begin{enumerate}
  \item Definable predicates~$\varphi$ formalize a notion of “explicit computability” of~$\varphi$, in a certain local and approximate sense.
    Namely, given (\emph{i})~any “approximation error” $\varepsilon>0$, and (\emph{ii})~some \emph{a priori} knowledge of the argument~$v$
    (i.e., knowing that $v$ belongs to a specific shard~$L[\rb]$ —this is the sense of~“locality” of the computation),
    one may regard the $\varepsilon$-uniformly approximating formula $\psi^{(\varepsilon)}_{[\rb]}$ to~$\varphi$ on~$L[\rb]$ as an explicit algorithm that (modulo an approximation error not exceeding~$\varepsilon$) computes~$\varphi(v)$.
    Numerical algorithms relying on arbitrary-precision operations are typically definable in the above sense:
    On the one hand, one must ensure that the calculation is stable under rounding errors (of the order of the machine’s~$\varepsilon$);
    on the other, such rounding errors on inputs may lead to arbitrarily large output error unless the magnitude of inputs is bounded
    (i.e., unless the inputs belong to a given shard)
    \emph{a priori}.

  \item By the definition of the topologies on~$L$ and $\LSh$, every monomial~$P$ is continuous and bounded by $r_P$ on any shard $L[\rb]$, and extends continuously to $\LSh$ (as the $P$-coordinate function).
    Since connectives are obtained by pointwise application of continuous real-valued functions of real arguments (addition and multiplication), every polynomial on~$L$ is also continuous, and extends to a continuous bounded function on type-shards~$\fL[\rb]$.
    Definable predicates, on the other hand, need not be continuous on~$L$
    —although their restrictions to shards $L[\rb]$ necessarily are continuous and bounded
    (being uniform limits of polynomials on~$L[\rb]$, which is compact).
  \item Let $\rb$ be an arbitrary sizer.
    The restriction of a monomial $P$ to the type-shard $\fL[\rb]$ admits the \emph{a priori} bound $C= r_P$, so that $P{\restriction}\fL[\rb]$ takes values in~$[-C,C] = [-r_P,r_P]$.%
    \footnote{A constant $r$ also admits the trivial bound $C = \nrm{r}$.}
    By recursion on the application of connectives leading to an arbitrary polynomial~$\varphi$ from monomials, \emph{a priori} bounds $C = C^{\varphi}_{\rb}\in[0,\infty)$ such that $\varphi{\restriction}\fL[\rb]$ takes values in $[-C,C]$ are easily found.
    (Recursively apply the rules:
    $C_{\rb}^{\varphi+\psi}\coloneqq C_{\rb}^{\varphi}+C_{\rb}^{\psi}$, and $C_{\rb}^{\varphi\psi}\coloneqq C_{\rb}^{\varphi}\cdot C_{\rb}^{\psi}$.)
  \item By definition of the topology on~$L$ and the Reduction Axioms, the collection of (continuous) predicates $P(\cdot): \LSh\to\RR$ (extended to the type space~$\LSh$) separates points of~$\LSh$
    (\emph{a fortiori}, points of any type-shard $\fL[\rb]$).
    By the Stone--Weierstrass Theorem, any Sh-extendable $\varphi: L\to\RR$ is necessarily definable.
    (In particular, any continuous $\varphi: \LSh\to\RR$ is definable in such case.)
    Clearly, the condition may be relaxed to requiring that $\varphi$ have continuous restrictions to type-shards $\fL[\rb]$ for $\rb$ in some exhaustive~$R$.
    By contrast, continuous predicates $L\to\RR$ need not be definable.
  \item In general, a function $\varphi:L\to\RR$ whose restrictions to shards are continuous need \emph{not} be continuous on~$L$
    (not even under the additional assumption that $L$ be realcompact).
    For $\cP$ (at most) countable, however, Sh-continuous real functions on the type space $\LSh$ are continuous
    (since $\LSh = \fL$ is a k-space in such case, by Proposition~\ref{prop:standard-types-countable-P}).
  \end{enumerate}
\end{remarks}

\subsubsection{Definable features}
\label{sec:definable-features}

Remarks~\ref{rem:disconts-ucomps} provide relevant context for this subsection.

Given $P\in\cP$, the \emph{$P$-feature} of a transition-in-type $\ff\in\fT$ is the real-valued function
\begin{align*}
  P\circ\ff: \LSh &\to \RR\\
  \fv &\mapsto P(\ff(\fv)).
\end{align*}
(One may call such a feature “atomic” or “monomial”.)

Individual features of a transition-in-type $\ff\in\fT$ may be definable or non-definable.
A transition-in-type is \emph{definable}
if its features are definable.

In the setting of §~\ref{sec:deep-equilibria}, one may ask under what circumstances a specific feature of a deep computation~$\ff\in\fD$ is effectively computable.

Sh-continuous features of transitions are definable, by Remark~\ref{rem:definability}(4).

\subsection{Definability of ultracomputations-in-type}
\label{sec:definable-ucomp}

Nonprincipal ultrafilters $\cU$ on infinite sets are ineffably inexplicit.
Thus, as a first step towards grasping ultracomputations, it is natural to consider ultralimits $\gammab^{\cU}$ of pointwise-bounded sequences $\gammab=(\gamma_n)_{n\in\NN}\subseteq\Gamma$ indexed by the infinite countable set~$\NN$.
Ultracomputations $\gammab^{\cU}$ obtained in this form (as $\gammab$ and~$\cU$ vary) are accumulation points of arbitrary \emph{countable} sets of realized computations.

Ultralimits obtained from countable subsets of~$\Gamma(\cdot)$, although less general than those obtained from arbitrary subsets, may still be quite complex.
Given a countable set $\gammab = (\gamma_n(\cdot))_{n<\omega}\subseteq\Gamma$ of pointwise-bounded computations, it is natural to consider \emph{sequential} limits of~$\gammab$, i.e., ultracomputations arising as pointwise limits of subsequences of $\gammab$, namely ultracomputations $\tilde{\gammab}$ of the form
\begin{equation*}
  \fv\mapsto \tilde{\gammab}(\fv) \coloneqq \lim_{k\to\infty}\gamma_{n_k}(\fv)
\end{equation*}
for subsequences (otherwise arbitrary) $\tilde{\gammab} = (\gamma_{n_k})_{k\in\NN}$ of~$\gammab$.

By Proposition~\ref{prop:sh2sh-compactness}, pointwise-boundedness of~$\gammab$ implies that \emph{all} ultracomputations $\gammab^{\cU}$ exist for arbitrary~$\cU$ on the index set~$I$ of any family $\gammab = (\gamma_i)_{i\in I}$
—regardless of the cardinality of~$\cP$ or~$L$.
Ultracomputations $\gammab^{\cU}$ realizable from \emph{sequences} $\gammab = (\gamma_n)_{n\in\NN}$ are quite special;
those realizable as \emph{sequential} limits~$\tilde{\gammab}$, even more so.

If $\gammab$ is a pointwise-bounded \emph{sequence}, \emph{and} $\cP$ is at most countable, then \emph{at every fixed $\fv\in \LSh$}, the ultralimit $\tilde{\gammab}(\fv)$ is realized as a sequential limit (by a standard diagonalization argument);
however, the realizing subsequence $(\gamma_{n_k})$ will typically depend on~$\fv$ and cannot be chosen uniformly over $\fv\in \LSh$.
When $\cP$ is uncountable, sequentially realizing an ultralimit of $\gammab$ —even at a single point~$\fv$— may be unfeasible.

The results in this concluding section %
relate (\emph{i})~continuity on shards of ultracomputations, (\emph{ii})~the ability to obtain such ultracomputations as accumulation points of \emph{countable sets} of computations, or as \emph{sequential limits} of computations, (\emph{iii})~the \emph{definability} of such ultracomputations, and~(\emph{iv}) a \emph{limit-exchange criterion} (originally due to Grothendieck).

\subsubsection{Relative compacta of continuous layer transitions}
\label{sec:relcpct-layer-trans}

For any topological space~$X$, let $\Cp(X)\subseteq\RR^X$ be the set of all continuous real functions on~$X$, endowed with the relative (subspace) topology of the product~$\RR^X$, i.e., the topology of point-wise convergence at each $x\in X$.
More generally, given two spaces~$X,Y$, the space $\Cp(X;Y)$ is the subspace of the product $Y^X = \prod_{x\in X}Y$ consisting of continuous maps $X\to Y$.
(“$\Cp$” means “pointwise topology on continuous functions”.)

Note that $\Cp(X), \Cp(X;Y)$ are generally \emph{not} closed subspaces of~$\RR^X, Y^X$.

A Hausdorff topological space~$Z$ is \emph{countably compact} if every infinite (equivalently, every infinite countable) subset~$B\subseteq Z$ has a limit point~$z\in Z$.
A subset $A\subseteq Z$ is \emph{relatively countably compact} (or \emph{countably compact in~$Z$}) if every infinite (equivalently, every infinite countable) subset~$B\subseteq A$ has a limit point~$z\in Z$.

(One may take the properties above as the definition of (relatively) countably compact for arbitrary, not necessarily Hausdorff spaces~$Z$.
However, the Hausdorff assumption implies desirable additional properties, e.g., \cite[Theorems~3.10.2, 3.10.3, etc.]{Engelking1989}.
In our applications, $Z$ is always a subspace of the layer state space~$L$ of a CCS, or of the type space~$\fL$, and hence Hausdorff.)

A topological space~$Y$ is \emph{angelic} if (\emph{i})~every relatively countably compact subset $A\subseteq Y$ is relatively compact, and (\emph{ii})~the closure $\overline{A}\subseteq Y$ of any such (relatively compact)~$A$ consists precisely of limits of sequences in~$A$.%
\footnote{A topological space possessing property (\emph{ii}) above is called \emph{Fréchet-Urysohn.}}

\subsubsection{A topological result of Grothendieck}
\label{sec:Grothendieck2}

\begin{therm}\label{thm:Grothendieck}
  Let
  \begin{itemize}
  \item $X$ be a countably compact topological space;
  \item $Y$, any Tychonoff space, having the property that its relatively countably compact subsets are relatively compact (which necessarily holds in case $Y$ is realcompact);
  \item $X_0\subseteq X$ any dense subset.
  \end{itemize}
    Then:
  \begin{enumerate}
  \item $\Cp(X;Y)$ is angelic.
  \item Assume that $Y$ is explicitly embedded as a subspace $Y\subseteq\RR^{\cP}$ for some index set~$\cP$.
 A set $F\subseteq\Cp(X;Y)$ of continuous maps $X\to Y$ is relatively compact if and only if
    \begin{enumerate}
    \item $F$ is pointwise bounded (i.e., $\{P(f(x)): f\in F\}$ is bounded for each $P\in \cP$ and $x\in X$)%
      \footnote{Here, we use the notation $P(f(x))$ for the “$P$-th coordinate” $f_P(x)$ of any $f\in(\RR^{\cP})^X$.}, and
    \item for all sequences $(f_m)_{m\in\NN}\subseteq F$, $(x_n)_{n\in\NN}\subseteq X_0$, any $P\in\cP$ and ultrafilters $\cU,\cV$ on~$\NN$, the following equality (called the \emph{limit-exchange property}) holds between iterated ultralimits:
      \begin{equation}\label{eq:limit-exchange-Groth}
        \Ulim_m\Vlim_nP(f_m(x_n))
        = \Vlim_n\Ulim_mP(f_m(x_n)),
      \end{equation}
      which both exist.
    \end{enumerate}

  \item Even if all hypotheses on~$X,Y$ pertaining to compactness are omitted (i.e., $Y$ is Tychonoff and $X$ arbitrary), the limit-exchange condition~(b) alone implies that every accumulation point of~$F\subseteq Y^X$ is continuous
    (i.e., the closure~$\overline{F}\subseteq \Cp(X;Y)$).
  \end{enumerate}
\end{therm}
  For a contemporary exposition of Grothendieck’s theorem and its consequences, we refer the reader to the paper on angelic spaces and the double limit relation by König and Kuhn~\cite{Konig-Kuhn:1987}.
\begin{proof}
  Theorem~\ref{thm:Grothendieck} aggregates several results in Grothendieck’s \emph{“Critères de compacité”} \cite[Théorèmes~1 \&~2, Remarque~2, Corollaire~2]{Grothendieck:1952}.
  Presently, we merely offer some remarks on translating between French terms and decades-old nomenclature to their contemporary equivalents in English.
  Spaces $\mathrm{C}_{\mathrm{s}}(X;Y)$ (where “s” refers to the “simple” topology, i.e., of pointwise convergence) are now denoted $\Cp(X;Y)$ (or just~$\Cp(X)$, when $Y=\RR$).
  \emph{“(Relativement) semi-compact”} (resp., \emph{“relativement compact”}) refers to \emph{(relatively) countably compact} (resp., relatively compact) sets.
  Functions take values in~$Y$, which we take to be a Tychonoff space (“complètement régulier”
  —i.e., completely regular \emph{and Hausdorff} in the standard contemporary sense) endowed with an embedding into a product $\RR^{\cP}$, hence $Y$ is a uniform Hausdorff space (\emph{“espace uniforme séparé”})~\cite[§§~1.5, 3.10, 8.1]{Engelking1989}.
\end{proof}

\begin{remarks}
  \begin{enumerate}
  \item Condition~(\emph{a}) above implies that both iterated ultralimits in equation~\eqref{eq:limit-exchange-Groth} in~(\emph{b}) exist.
    However, (\emph{b}) explicitly \emph{asserts} the requirement that the limits exist
    —not merely that they are equal \emph{when} they exist.
  
  \item The hypotheses on~$Y$ are satisfied if $Y$ is realcompact, in which case the embedding $Y\subseteq\RR^{\cP}$ is as a closed subspace of the product;
   moreover, any $Y$ embedded as a closed subspace of any such product of lines satisfies all hypotheses (including those in part~(2) of the theorem).
  \end{enumerate}
\end{remarks}

\subsubsection{The Fundamental Theorem of Definability}
\label{sec:fundamental-thm}

\begin{therm}\label{thm:fundl-thm-defin}
  Let $\cC = \langle \underline{L},\underline{\Gamma}\rangle$ be a CCS.
  Let $R$ be an exhaustive sizer collection, and let $\Delta\subseteq\Gamma[R]$ be any $R$-confined set (of Sh-extendable computations, by assumption).
  Then, the properties below are equivalent:

\paragraph{\textbf{Extendable Ultracomputations (uExt)}}
Every ultracomputation over~$\Delta$ is Sh-extendable.

\paragraph{\textbf{Limit Exchange (LE)}}
 For all sizers $\rb$, all sequences $v_{\bullet}\subseteq L[\rb]$ and $\gammab\subseteq \Delta$, and ultrafilters $\cU,\cV$ on~$\NN$, the iterated ultralimits $\Ulim_m\Vlim_n \gamma_m(v_n)$ and $\Vlim_n\Ulim_m \gamma_m(v_n)$ both exist and are equal:
\begin{equation}\label{eq:limit-exchange2}
  \Ulim_m\Vlim_n \gamma_m(v_n) = 
  \Vlim_n\Ulim_m \gamma_m(v_n).
\end{equation}

\paragraph{\textbf{Uniform Approximation (UA)}}
Every ultracomputation $\ff$ over~$\Delta$ is definable without parameters:
For any sizer $\rb$, any $\varepsilon>0$, and all $P\in\cP$, there exists a polynomial~$\psi = \psi_{\rb,\varepsilon,P}$ (without parameters) such that
\begin{equation}\label{eq:unif-approxn}
  \nrm{\psi(v) - P(\ff(v))} < \varepsilon\qquad
  \text{for all $v\in L[\rb]$.}
\end{equation}
Moreover:
\begin{enumerate}
\item In case any (hence all) of the above conditions hold for~$\Delta$,
  the restriction of any ultracomputation $\ff$ over~$\Delta$ to any type-shard~$\fL[\rb]$ is the limit $\ff{\restriction}\fL[\rb] = \lim_n\gamma_n(\cdot){\restriction}\fL[\rb]$ obtained from a sequence $\gammab \subseteq \Delta$.
\item For arbitrary $\Delta\subseteq\Gamma$
  (i.e., $\Delta$ not \emph{a priori} included in~$\Gamma[R]$ for some exhaustive~$R$), the Limit Exchange condition alone implies that all ultracomputations
  over~$\Delta$ are Sh-extendable.%
    \footnote{The explicit LE hypothesis that both iterated ultralimits in~\eqref{eq:limit-exchange2} exist is essential when $R$ and the implied pointwise bounds on~$\Delta$ are not given \emph{a priori}.}
\end{enumerate}
\end{therm}
\begin{proof}[Proof of Theorem~\ref{thm:fundl-thm-defin}]
  Because of the hypothesis $\Delta\subseteq\Gamma[R]$, it is quite clear that one may specialize all uses of sizers~$\rb$ and universal properties of sizers to involve sizers $\rb\in R$ only.

  In Grothendieck’s Theorem~\ref{thm:Grothendieck}, let $Y = \fL \subseteq \RR^{\cP}$ (realcompact) and, for a momentarily fixed $\rb\in R$, let $X = \fL[\rb]$ (compact, hence countably compact),
  and $Z \coloneqq \Cp(\fL[\rb];\fL)$.
  Denote by $\Delta[\rb] \subseteq Z$ the set of functions~$\gamma_{[\rb]}\coloneqq \gamma(\cdot){\restriction}\fL[\rb]$ as $\gamma\in\Delta$ varies.
  By Theorem~\ref{thm:Grothendieck}, %
  the condition that all ultracomputations over~$\Delta$ are continuous on $\fL[\rb]$ is equivalent to the relative compactness of $\Delta[\rb]\subseteq Z$.

  Since $Z$ is angelic (Theorem~\ref{thm:Grothendieck}(1)), assertion (1) follows.

  The pointwise boundedness condition~2(a) in Theorem~\ref{thm:Grothendieck} is satisfied since $\Delta[\rb]$ is pointwise bounded (as $\Delta$ is uniformly confined by assumption);
  therefore, relative compactness of~$\Delta[\rb]$ is, in turn, characterized by the Limit Exchange condition (equivalent to 2(b)), so LE is equivalent to the preceding three conditions.
  Moreover, assertion (2) follows from Theorem~\ref{thm:Grothendieck}(3).

  Any feature $P\circ\ff: \LSh\to\RR$ of any transition-in-type~$\ff:\LSh\to\LSh$, if uniformly approximable on some shard~$L[\rb]$ by polynomials~$\psi$
  —any of which has a unique extension to a continuous real function on~$\fL$, bounded on~$\fL[\rb]$—
  must necessarily extend continuously to~$\fL[\rb]$.
  Letting $P\in\cP$ and $\rb$ vary, we see that a definable ultracomputation is necessarily Sh-continuous:
  UA implies uExt.
  Reciprocally, by the Stone--Weierstrass Theorem, every continuous real function $\fL[\rb]\to\RR$ is uniformly approximable by polynomials in predicates $P\in\cP$ (because these predicates separate points of~$\fL[\rb]$), i.e., by polynomials without parameters.
  Therefore, any Sh-continuous ultracomputation is definable without parameters:
  uExt implies UA.  
\end{proof}

\begin{remarks}\label{rem:Grothendieck}
  \begin{enumerate}
  \item The extendability condition (uExt) in Theorem~\ref{thm:fundl-thm-defin} may be regarded as auxiliary in proving the equivalence LE$\Leftrightarrow$UA.
    The implication UA$\Rightarrow$LE is not difficult to prove directly:
    On the one hand, UA$\Rightarrow$uExt by the straightforward argument in the proof above.
    Afterward, uExt$\Rightarrow$LE follows easily:
    uExt implies that every ultracomputation $\ff\coloneqq\Ulim_n\gamma_n(\cdot)$ is continuous on any compact~$\fL[\rb]$, and LE simply states the continuity of $\ff$ at ultralimit points of the form $\mathfrak{v}\coloneqq\Vlim_n\tp(v_n)\in\fL[\rb]$ for arbitrary state sequences $(v_n)\subseteq L[\rb]$.

    By contrast, the implication LE$\Rightarrow$UA may be seen as a significantly deeper consequence of Grothendieck’s Theorem:
    A natural limit-exchange condition implies that layer transitions-in-type are explicitly computable!

  \item One could take a probabilistic approach to the uniqueness and computability of equilibria inspired by ideas from deep learning and the Examples~\ref{exa:DEQ-unit-interval} and~\ref{exa:DEQ-Newton} in §~\ref{sec:deep-equi-exa}.
    For simplicity, assume that $L$ is realcompact (so $L = \LSh = \fL$).
    The uniqueness and continuity of deep iterates $\gamma^{*}$ at a state~$v\in L$ may be tested empirically by taking finitely many independent random points $(v_i)_{i<k}$ in a small neighborhood of~$v$ and computing $w_i = \gamma^{n_i}(v_i)$ for some large and also random integers $(n_i)_{i<k}$.
  To the extent that the points $(w_i)_{i<k}$ are (or are not) near each other, one may infer (in a statistical sense) whether $\gamma^{*}$ is (or is not) continuous at~$v$ with increasingly larger probability as $k$ grows.
  At points of continuity $v$ (as determined with high probability taking $k$ sufficiently large), any of the computed points $w_i$ may be regarded as an approximation to the exact and unique value $\gamma^{*}(v)$.
  This approach hints at a relativized notion of computability based on almost-everywhere (or at least local) continuity rather than everywhere continuity, which we intend to revisit in a sequel paper.
  \end{enumerate}
\end{remarks}

\appendix

\section{Smooth ultracomputations and effectively computable equilibria in deep neural networks}
\label{sec:def-equi}

Extending the framework of the main body of the paper, one may introduce \emph{smooth (ultra)computations} as those having output features varying smoothly (i.e., differentiably) with the input features.
Considerations of differentiability —particularly in infinite dimensions— are very delicate and exceed the scope of this current paper
(after all, our notions of extendability and definability only capture continuity properties).
Since differentiability is an essential assumption in current approaches to effective/implicit computability of deep neural networks, this appendix is a brief and informal outline of extensions to our framework beyond the present topological context so as to capture differentiability.

\emph{Throughout this appendix, we fix a \emph{realcompact} CCS $\cC$ whose layer state space $\underline{L}$ is a differentiable (smooth) manifold of finite dimension~$n$, and all predicates $P\in\cP$ are differentiable on~$L$.}

In particular, we assume that the embedding $L\subseteq\RR^{\cP}$ is as a closed subspace (in the product topology).

\subsection{Deep equilibria of neural networks \emph{à la} Bai--Kolter--Koltun}
\label{sec:BaiEtAl}

\subsubsection{Unique Deep Equilibria}
An empirical observation in the context of Neural Network Deep Equilibrium Models~\cite{BaiKolterKoltun2020} is that, in situations where a deep iterate $\gamma^{(\cU)} = \Ulim_n\gamma^n$ of some computation~$\gamma$ (assumed confined, for simplicity) exists, it is often \emph{independent} of the ultrafilter~$\cU$.%
\footnote{Implicitly, both \cite{Chen-Rubanova-Bettencourt-Duvenaud:2018} and \cite{BaiKolterKoltun2020} work in a setting where the states space $L=\fL$ is realcompact, so there is no distinction between transforms and transitions-in-type.}
In such case, all deep iterates $\gamma^{(\cU)}$ are one and the same transition $\gamma^{*}: L \to L$
—a “deep state” of the NN obtained by iteration of~$\gamma$.
Therefore, the sequence of iterates $\gamma^n$ converges pointwise to the t-t~$\gamma^{*}$ as an ordinary limit (rather than only as an ultralimit).
We say that such $\gamma$ has the \emph{Unique Deep Equilibrium (UDEQ) Property}. Smoothness properties of $\gamma$ are required for important applications, as described below. 

\subsubsection{Fixed-point algorithms as “Black Boxes”}
\label{sec:black-box}

Bai \emph{et al.}\ note (empirically) that, for NNs obtained by iterating a common “weight-tied” layer transition~$\gamma$, the deep state $\gamma^{*}$ takes any input state $v\in L$ to another $v^{*} = \gamma^{*}(v)$ that is fixed by (the t-t implied by) the original~$\gamma$, i.e., $\gamma(v^{*}) = v^{*}$;
in other words, $\gamma^{*}$ takes values in the set $\Fix(\gamma) = \{v\in L: \gamma(v) = v\}$, so $\gamma^{*}$ is a deep equilibrium (DEQ) in a very strong sense.
Empirical findings also suggest that, given $\gamma\in\Gamma$, the DEQ state $\gamma^{*}: L\to \Fix(\gamma)\subseteq L$ %
may be well approximated by some generic “black-box” fixed-point algorithm \FFix.
Such an algorithm should take as inputs the transformation $\gamma$ and initial state-in-type~$v$, and return the fixed point $v^{*}=\gamma^{*}(v) = \FFix(\gamma;v)$.

Like any algorithm based on arbitrary-precision arithmetic, what such an algorithm $\FFix$ does in practice, given an acceptable error $\varepsilon>0$ and finitely many output features $Q_1,\dots,Q_k\in\cP$ specified in advance, is to return a suitable $k$-tuple $\FFix_k(\gamma;v; \varepsilon) = (r_1,\dots,r_k)$ of real numbers such that, for $1\le i\le k$, $\nrm{Q_i(\gamma^{*}(v))-r_i}<\varepsilon$.
Under our current assumption that $L$ is smooth of dimension~$n$, \emph{all} features $P\in\cP$ of the output are (heuristically speaking, and perhaps only locally) implicitly defined in terms of some $n$-many input features $Q_i$, $1\le i\le n$.
Moreover, a generic such algorithm \FFix\ typically assumes that the given map $\gamma$ is not merely continuous but smooth (or at least sufficiently differentiable), and relies on gradient-based methods.

In principle, evaluating (or, at least, approximating) the map $v\mapsto\gamma^{*}(v)$ by means of a “black-box” \FFix\ results in comparable computational complexity or even savings over the direct method of computing successive iterates $\gamma^n(v)$ until a limit is (very nearly) reached.
Memory savings in training DEQ networks (cf.\ §~\ref{sec:PFC} below) are also a key advantage to their success.
From a theoretical perspective, the innovation lies in effectively bringing deep networks (at least, when obtainable as iterative deep equilibria) on par with classical networks, thereby enriching the class of directly and efficiently computable functions.

\subsubsection{Parametrized Families of Computations}
\label{sec:PFC}

Fix a CCS $\underline{\cC} = \langle \underline{L}, \underline{\Gamma} \rangle$ with underlying CSS $\underline{L} = \langle L,\cP\rangle$ as layer state space, as well as a second CSS $\underline{X} = \langle X, \cQ\rangle$, called the space of \emph{computation parameters}, and a map $F: X \to \Gamma: x\mapsto F_x$, which we regard as a parametrization of (some) computations by elements (parameters) $x\in X$.
We make the same assumptions about~$\underline{X}$ as about~$\underline{L}$ above (namely, $\underline{X}$ is a finite-dimensional differentiable manifold embedded as a closed subspace of~$\RR^{\cQ}$).
We call the structure $\underline{F} = \langle F, \underline{L}, \underline{X}, \underline{\Gamma}\rangle$ a \emph{Parametrized Family of Computations (PFC)} (all of which are confined).
It is appropriate to think of the parameter $x\in X$ as the “weights” of the computation~$F_x$.

We assume that $\Gamma$ has only confined transitions.
It is quite natural to assume that $F$ is (\emph{i})~continuous (as a map $X\to\fT$), and (\emph{ii})~confined, i.e., that it maps each shard of the parameter type space~$\fX$ into $\fT[\sbb^{(\cdot)}]$ for some collection~$\sbb^{(\cdot)}$ of pointwise bounds.

A UDEQ hypothesis for $\underline{F}$ implies a map  $X\to\fT: x\mapsto F^{*}_x$ which may also be regarded as a map
\begin{align*}
  F^{*}: X\times L &\to L\\
  (x,v) &\mapsto F^{*}(x;v).
\end{align*}

\subsubsection{Training deep networks}
\label{sec:training}

Training the deep neural network $F^{*}_x$ translates to finding weights $x$ such that $F^{*}_x$ satisfies a given condition, which we presently take to mean minimizing a given/specified real-valued loss function $\ell:\fT\to[0,\infty)$.
(At least intuitively, if not necessarily literally, the value $\ell(\fg)\ge0$ captures how far a transition $\fg\in\fT$ is from an optimal/idealized $G\in\fT$.)
Since $F^{*}_x$ for fixed $x$ is implicitly defined by either a fixed-point condition or ODE as above, the enormous memory cost of back-propagation through layers%
\footnote{Not least, because back-propagation would involve an unbounded number of ordinary layers to begin with!}
is replaced by that of minimizing the function $\tilde{\ell} \coloneqq \ell\circ F^{*}: X\to[0,\infty)$.
Note that $\tilde{\ell}$ is merely a new real-valued predicate on the parameters CSS $\underline{X}$.
Assuming that $\tilde{\ell}$ is shard-continuous, it is definable, hence depends \emph{de facto} on only finitely many features $Q_1,\dots,Q_n\in\cQ$ of its input $x\in X$ (up to an arbitrarily small admissible error $\varepsilon>0$).
Assuming $\tilde{\ell}$ is smooth as well, the deep network may be trained using standard/“black-box” gradient-based procedures to find a minimizer $x\in X$ for $\tilde{\ell}$.
However, we note that it is essential for $\tilde{\ell}$ to be definable in order to allow even the possibility that some algorithm involving arbitrary-precision arithmetic and finitely many real quantities at a time succeeds in finding the minimizer.

\subsection{Neural ODEs \emph{à la} Chen--Rubanova--Bettencourt--Duvenaud}
\label{sec:ChenEtAl}

In another setting that is technically different but conceptually closely related to the one in~§~\ref{sec:BaiEtAl}, Chen~\emph{et al.}~\cite{Chen-Rubanova-Bettencourt-Duvenaud:2018} also model deep states of residual networks (“Neural ODEs”) using differential equation techniques.
The intuition behind Neural ODEs is the following:
Consider a layered computation with atomic-step sequence $\gammab = (\gamma_0,\gamma_1,\dots)$ such that all such steps $\gamma_i$ are “residually” very small (in the sense that the input and output features of any atomic step~$\gamma_i$ differ very little).
Successive $n$-composites $\gammab^{(n)} = \gamma_{n-1}\dots\gamma_1\gamma_0$ change very little with~$n$;
as one varies $\gammab$ in such a way that the atomic steps residually vanish (i.e., $\gamma_i$ is vanishingly close to~$\id$) and one allows $n$ to grow without bound, when a limit exists, Chen \emph{et al.}\ model it as a family $(\gamma^{(t)})_{t\ge 0}$ (indexed by a real “time” variable~$t\ge0$) of transitions $\gamma^{(t)}$, which we assume to be (confined) elements of~$\fT$.
(The real variable $t$ captures an appropriate asymptotic rescaling of the “discrete time”~$n$.)
In this manner, each value $t=t_0$ captures a specific notion of deep state (as an asymptotic limit of deep composites of residually small layered transitions), realized as a confined transition.

One may hope that such transitions $\gamma^{(t)}$ vary differentiably with~$t$;
this suggests modeling the entire family $(\gamma^{(t)})_{t\ge 0}$ of deep computational states per the differential equation implied.
(In this manner, for each fixed $t=t_0\ge0$, one obtains a deep network $\gamma^{(t_0)}$ in some sense.)

Thus, “Neural ODEs” arise from differential equations of the form
\begin{equation}
  \label{eq:NeuralODE}
  \dot{v} = \mathbf{s}(v;t),
\end{equation}
where $\mathbf{s}: L\times[0,+\infty)\to \mathrm{T}L$ is a (time-dependent) section of the tangent bundle $\mathrm{T}L$ of the layer state space~$L$ (i.e., $\mathbf{s}(v;t)\in \mathrm{T}_vL$ for all $v\in L$ and $t\ge0$, where $\mathrm{T}_vL$ is the tangent space of~$L$ at~$v$). Interpreting the ODE~\eqref{eq:NeuralODE} hinges on the smooth manifold structure assumed of the state space~$L$.%
\footnote{The notions of differentiable structure and tangent space on an arbitrary layer space $L$ are neither well nor uniquely defined when $L$ is not finite-dimensional;
their formalization would require much stronger assumptions on~$L$, as well as the formalism of Banach spaces for tangent spaces~$\mathrm{T}_vL$.}
Chen \emph{et al.}~illustrate empirically the feasibility and effectiveness of modeling deep equilibria by Neural ODEs.
Let us denote the time-$t$ evolution by~\eqref{eq:NeuralODE} using the (hopefully, suggestive) notation $v\mapsto e^{t\mathbf{s}}(v)$, i.e., $e^{t\mathbf{s}}$ is the deep equilibrium of the Neural ODE solving~\eqref{eq:NeuralODE} (i.e., “$\gamma^{(t)}$” in the earlier informal discussion).%
\footnote{When the section $\mathbf{s} = \mathbf{s}(v)$ depends only on the state~$v$ (not on time~$t$), the ODE~\eqref{eq:NeuralODE} is autonomous.
A time evolution $e^{t \mathbf{s}}$ of such an autonomous Neural ODE is analogous to a “parameter-tied” deep equilibrium after Bai \emph{et al.}}
Effective computation of $e^{t\mathbf{s}}$ relies on a generic “black-box” ODE solver algorithm~\ODEsolve.
Such an algorithm should take as inputs the section $\mathbf{s}$, initial state-in-type~$v$ and time $t\ge0$, and return the output $e^{t \mathbf{s}}(v) = \ODEsolve(\mathbf{s};v,t)$.
(More realistically, such an \ODEsolve\ algorithm presumably would return approximate values for any finitely many specified features of~$e^{t \mathbf{s}}(v)$;
refer to the discussion of $\FFix_k$ above.)

Modeling deep computations by Neural ODEs and realizing them by means of an ODE solver effectively brings them on computational par with classical neural networks.
The key insight of Chen \emph{et al.}\ (which predates the work of Bai \emph{et al.})\ is that training such Neural ODEs may be done using the “adjoint sensitivity” method of Pontryagin instead of doing (extremely memory-intensive) back-propagation through layers
—which, at any rate, have been essentially abstracted away.
The adjoint sensitivity method may be implemented using \ODEsolve\ itself, so the training is both memory- and computation-efficient.
Formalizing their method to train Neural ODEs in the spirit of~§~\ref{sec:PFC} above requires parametrizing sections~$\mathbf{s}$ by a second CSS~$\underline{X}$; we omit the details.

\bibliographystyle{msclike}
\bibliography{biblioteca, iovino}


\end{document}